\documentclass[12pt]{article}
\usepackage{theorem}
\usepackage{amssymb}
\usepackage{graphicx}
\usepackage[mathscr]{eucal}
\newtheorem{prop}{}[section]

{\theorembodyfont{\upshape} \newtheorem{rema}[prop]{}}
\newcommand{\boma}[1]{{\mbox{\boldmath $#1$} }}

\begin{document}
\newcommand{\binom}[2]{\left( \barray{c} #1 \\ #2 \farray \right)}
\newcommand{\uper}[1]{\stackrel{\barray{c} {~} \\ \mbox{\footnotesize{#1}}\farray}{\longrightarrow} }
\newcommand{\nop}[1]{ \|#1\|_{\piu} }
\newcommand{\no}[1]{ \|#1\| }
\newcommand{\nom}[1]{ \|#1\|_{\meno} }
\newcommand{\UU}[1]{e^{#1 \AA}}
\newcommand{\UD}[1]{e^{#1 \Delta}}
\newcommand{\bb}[1]{\mathbb{{#1}}}
\newcommand{\HO}[1]{\bb{H}^{{#1}}}
\newcommand{\Hz}[1]{\bb{H}^{{#1}}_{\zz}}
\newcommand{\Hs}[1]{\bb{H}^{{#1}}_{\ss}}
\newcommand{\Hg}[1]{\bb{H}^{{#1}}_{\gam}}
\newcommand{\HM}[1]{\bb{H}^{{#1}}_{\so}}
\newcommand{\vers}[1]{\widehat{#1}}
\def\tvainf{\vspace{-0.4cm} \barray{ccc} \vspace{-0,1cm}{~}
\\ \vspace{-0.2cm} \longrightarrow \\ \vspace{-0.2cm} \scriptstyle{T \vain + \infty} \farray}
\def\k{t}
\def\iep{(-1,+\infty)}
\def\ep{\xi}
\def\Do{\mathscr{E}}
\def\lam{\lambda}
\def\Lam{\Lambda}
\def\vh{\vers{h}}
\def\vk{\vers{k}}
\def\er{\epsilon}
\def\erd{\er_0}
\def\op{\,\mbox{\scriptsize{or}}\,}
\def\S{\mathscr{H}}
\def\Cgot{{\mathfrak C}}
\def\Mn{\Gamma^{<}_n}
\def\Mpn{\Gamma^{\geqs}_n}
\def\pp{p'}
\def\ps{p''}
\def\Pip{P'}
\def\Pis{P''}
\def\wp{w'}
\def\ws{w''}
\def\Wp{W'}
\def\Ws{W''}
\def\Ed{\hat{E}}
\def\Dd{\hat{D}}
\def\um{u_{-}}
\def\up{\mathfrak{u}}
\def\el{t}
\def\em{z}
\def\uu{\lambda}
\def\dK{\delta {\mathscr K}}
\def\dG{\delta {\mathscr G}}
\def\DK{\Delta {\mathscr K}}
\def\DG{\Delta {\mathscr G}}
\def\Ti{{\mathscr T}}
\def\Km{{\mathscr K}}
\def\Ll{\mathscr{L}}
\def\Hh{\mathscr{H}}
\def\Mm{{\mathscr M}}
\def\Nn{{\mathscr N}}
\def\Rr{{\mathscr R}}
\def\RRr{{\mathcal R}}
\def\Gg{{\mathscr G}}
\def\Zz{Z}
\def\Ss{{\mathscr S}}
\def\Fe{{\mathscr F}}
\def\Ei{{\mathscr E}}
\def\Ww{{\mathscr W}}
\def\we{\wedge}
\def\We{\bigwedge}
\def\dbar{\hat{d}}
\def\Cc{\mathscr{C}}
\def\SZ{\mathcal{S}}
\def\TZ{\mathfrak{S}}
\def\CQ{S}
\def\GQ{\mathfrak{G}}
\def\C{{C_d}}
\def\Ac{\overline{A}}
\def\Bc{\overline{B}}
\def\Xd{X^{- \, d}_k}
\def\Yd{X^{+ \, d}_k}
\def\comple{\scriptscriptstyle{\complessi}}
\def\nume{0.407}
\def\numerob{0.00724}
\def\deln{7/10}
\def\delnn{\dd{7 \over 10}}
\def\e{c}
\def\p{p}
\def\z{z}
\def\symd{{\mathfrak S}_d}
\def\del{\omega}
\def\Del{\delta}
\def\Di{\Delta}
\def\mmu{\hat{\mu}}
\def\rot{\mbox{rot}\,}
\def\curl{\mbox{curl}\,}
\def\XS{\boma{x}}
\def\TS{\boma{t}}
\def\DS{\boma{\rho}}
\def\KS{\boma{k}}
\def\LS{\boma{\lambda}}
\def\PR{\boma{p}}
\def\VS{\boma{v}}
\def\ski{\! \! \! \! \! \! \! \! \! \! \! \! \! \!}
\def\h{L}
\def\EM{M}
\def\EMP{M'}
\def\R{R}
\def\E{E}
\def\FFf{\mathscr{F}}
\def\A{F}
\def\Xim{\Xi_{\meno}}
\def\Ximn{\Xi_{n-1}}
\def\lan{\lambda}
\def\om{\omega}
\def\Om{\Omega}
\def\Sim{\Sigm}
\def\Sip{\Delta \Sigm}
\def\Sigm{{\mathscr{S}}}
\def\Ki{{\mathscr{K}}}
\def\Hi{{\mathscr{H}}}
\def\zz{{\scriptscriptstyle{0}}}
\def\ss{{\scriptscriptstyle{\Sigma}}}
\def\gam{{\scriptscriptstyle{\Gamma}}}
\def\so{\ss \zz}
\def\Dz{\bb{\DD}'_{\zz}}
\def\Ds{\bb{\DD}'_{\ss}}
\def\Dsz{\bb{\DD}'_{\so}}
\def\Dg{\bb{\DD}'_{\gam}}
\def\Ls{\bb{L}^2_{\ss}}
\def\Lg{\bb{L}^2_{\gam}}
\def\bF{{\bb{V}}}
\def\Fz{\bF_{\zz}}
\def\Fs{\bF_\ss}
\def\Fg{\bF_\gam}
\def\Pre{P}
\def\UUU{{\mathcal U}}
\def\fiapp{\phi}
\def\PU{P1}
\def\PD{P2}
\def\PT{P3}
\def\PQ{P4}
\def\PC{P5}
\def\PS{P6}
\def\Q{P6}
\def\X{Q2}
\def\Xp{Q3}
\def\Vi{V}
\def\bVi{\bb{V}}
\def\K{V}
\def\Ks{\bb{\K}_\ss}
\def\Kz{\bb{\K}_0}
\def\KM{\bb{\K}_{\, \so}}
\def\HGG{\bb{H}^\G}
\def\HG{\bb{H}^\G_{\so}}
\def\EG{{\mathfrak{P}}^{\G}}
\def\G{G}
\def\de{\delta}
\def\esp{\sigma}
\def\dd{\displaystyle}
\def\LP{\mathfrak{L}}
\def\dive{\mbox{div}\,}
\def\la{\langle}
\def\ra{\rangle}
\def\um{u_{\meno}}
\def\uv{\mu_{\meno}}
\def\Fp{ {\textbf F_{\piu}} }
\def\Ff{ {\textbf F} }
\def\Fm{ {\textbf F_{\meno}} }
\def\piu{\scriptscriptstyle{+}}
\def\meno{\scriptscriptstyle{-}}
\def\omeno{\scriptscriptstyle{\ominus}}
\def\Tt{ {\mathscr T} }
\def\Xx{ {\textbf X} }
\def\Yy{ {\textbf Y} }
\def\VP{{\mbox{\tt VP}}}
\def\CP{{\mbox{\tt CP}}}
\def\cp{$\CP(f_0, t_0)\,$}
\def\cop{$\CP(f_0)\,$}
\def\copn{$\CP_n(f_0)\,$}
\def\vp{$\VP(f_0, t_0)\,$}
\def\vop{$\VP(f_0)\,$}
\def\vopn{$\VP_n(f_0)\,$}
\def\vopdue{$\VP_2(f_0)\,$}
\def\leqs{\leqslant}
\def\geqs{\geqslant}
\def\mat{{\frak g}}
\def\tG{t_{\scriptscriptstyle{G}}}
\def\tN{t_{\scriptscriptstyle{N}}}
\def\TK{t_{\scriptscriptstyle{K}}}
\def\CK{C_{\scriptscriptstyle{K}}}
\def\CN{C_{\scriptscriptstyle{N}}}
\def\CG{C_{\scriptscriptstyle{G}}}
\def\CCG{{\mathscr{C}}_{\scriptscriptstyle{G}}}
\def\tf{{\tt f}}
\def\ti{{\tt t}}
\def\ta{{\tt a}}
\def\tc{{\tt c}}
\def\tF{{\tt R}}
\def\P{{\mathscr P}}
\def\V{{\mathscr V}}
\def\TI{\tilde{I}}
\def\TJ{\tilde{J}}
\def\Lin{\mbox{Lin}}
\def\Hinfc{ H^{\infty}(\reali^d, \complessi) }
\def\Hnc{ H^{n}(\reali^d, \complessi) }
\def\Hmc{ H^{m}(\reali^d, \complessi) }
\def\Hac{ H^{a}(\reali^d, \complessi) }
\def\Dc{\DD(\reali^d, \complessi)}
\def\Dpc{\DD'(\reali^d, \complessi)}
\def\Sc{\SS(\reali^d, \complessi)}
\def\Spc{\SS'(\reali^d, \complessi)}
\def\Ldc{L^{2}(\reali^d, \complessi)}
\def\Lpc{L^{p}(\reali^d, \complessi)}
\def\Lqc{L^{q}(\reali^d, \complessi)}
\def\Lrc{L^{r}(\reali^d, \complessi)}
\def\Hinfr{ H^{\infty}(\reali^d, \reali) }
\def\Hnr{ H^{n}(\reali^d, \reali) }
\def\Hmr{ H^{m}(\reali^d, \reali) }
\def\Har{ H^{a}(\reali^d, \reali) }
\def\Dr{\DD(\reali^d, \reali)}
\def\Dpr{\DD'(\reali^d, \reali)}
\def\Sr{\SS(\reali^d, \reali)}
\def\Spr{\SS'(\reali^d, \reali)}
\def\Ldr{L^{2}(\reali^d, \reali)}
\def\Hinfk{ H^{\infty}(\reali^d, \KKK) }
\def\Hnk{ H^{n}(\reali^d, \KKK) }
\def\Hmk{ H^{m}(\reali^d, \KKK) }
\def\Hak{ H^{a}(\reali^d, \KKK) }
\def\Dk{\DD(\reali^d, \KKK)}
\def\Dpk{\DD'(\reali^d, \KKK)}
\def\Sk{\SS(\reali^d, \KKK)}
\def\Spk{\SS'(\reali^d, \KKK)}
\def\Ldk{L^{2}(\reali^d, \KKK)}
\def\Knb{K^{best}_n}
\def\sc{\cdot}
\def\x{\mbox{{\tt x}}}
\def\g{ {\textbf g} }
\def\QQQ{ {\textbf Q} }
\def\AAA{ {\textbf A} }
\def\gr{\mbox{gr}}
\def\sgr{\mbox{sgr}}
\def\loc{\mbox{loc}}
\def\PZ{{\Lambda}}
\def\PZAL{\mbox{P}^{0}_\alpha}
\def\epsilona{\epsilon^{\scriptscriptstyle{<}}}
\def\epsilonb{\epsilon^{\scriptscriptstyle{>}}}
\def\lgraffa{ \mbox{\Large $\{$ } \hskip -0.2cm}
\def\rgraffa{ \mbox{\Large $\}$ } }
\def\restriction{\upharpoonright}
\def\m{m}
\def\Fre{Fr\'echet~}
\def\I{{\mathcal N}}
\def\ap{{\scriptscriptstyle{ap}}}
\def\fiap{\varphi_{\ap}}
\def\dfiap{{\dot \varphi}_{\ap}}
\def\DDD{ {\mathfrak D} }
\def\BBB{ {\textbf B} }
\def\EEE{ {\textbf E} }
\def\GGG{ {\textbf G} }
\def\TTT{ {\textbf T} }
\def\KKK{ {\textbf K} }
\def\HHH{ {\textbf K} }
\def\FFi{ {\bf \Phi} }
\def\GGam{ {\bf \Gamma} }
\def\sc{ {\scriptstyle{\bullet} }}
\def\a{a}
\def\c{\kappa}
\def\parn{\par\noindent}
\def\teta{M}
\def\elle{L}
\def\ro{\rho}
\def\al{\alpha}
\def\alc{\overline{\al}}
\def\dal{\mathfrak{a}}
\def\si{\sigma}
\def\be{\beta}
\def\dbe{\mathfrak{b}}
\def\bec{\overline{\be}}
\def\dbec{\overline{\dbe}}
\def\ga{\gamma}
\def\te{\vartheta}
\def\tet{\vartheta}
\def\teta{\theta}
\def\ch{\chi}
\def\et{\eta}
\def\complessi{{\bf C}}
\def\len{{\bf L}}
\def\reali{{\bf R}}
\def\interi{{\bf Z}}
\def\Z{{\bf Z}}
\def\naturali{{\bf N}}
\def\To{ {\bf T} }
\def\Td{ {\To}^d }
\def\Rt{ \reali^3 }
\def\Tt{ {\To}^3 }
\def\Zd{ \interi^d }
\def\Zt{ \interi^3 }
\def\Zet{{\mathscr{Z}}}
\def\Ze{\Zet^d}
\def\T1{{\textbf To}^{1}}
\def\Sfe{ {\bf S} }
\def\Sd{\Sfe^{d-1}}
\def\St{\Sfe^{2}}
\def\es{s}
\def\FF{\mathcal F}
\def\FFu{ {\textbf F_{1}} }
\def\FFd{ {\textbf F_{2}} }
\def\GG{{\mathcal G} }
\def\EE{{\mathcal E}}
\def\KK{{\mathcal K}}
\def\PP{{\mathcal P}}
\def\Pp{{\mathcal P}'}
\def\Ps{{\mathcal P}''}
\def\PPP{{\mathscr P}}
\def\PN{{\mathcal P}}
\def\PPN{{\mathscr P}}
\def\QQ{{\mathcal Q}}
\def\J{J}
\def\Np{{\hat{N}}}
\def\Lp{{\hat{L}}}
\def\Jp{{\hat{J}}}
\def\Vp{{\hat{V}}}
\def\Ep{{\hat{E}}}
\def\Gp{{\hat{G}}}
\def\Kp{{\hat{K}}}
\def\Ip{{\hat{I}}}
\def\Tp{{\hat{T}}}
\def\La{\Lambda}
\def\Ga{\Gamma}
\def\Si{\Sigma}
\def\Upsi{\Upsilon}
\def\Gam{\Gamma}
\def\Gag{{\check{\Gamma}}}
\def\Lap{{\hat{\Lambda}}}
\def\Upsig{{\check{\Upsilon}}}
\def\Kg{{\check{K}}}
\def\ellp{{\hat{\ell}}}
\def\j{j}
\def\jp{{\hat{j}}}
\def\BB{{\mathcal B}}
\def\LL{{\mathcal L}}
\def\MM{{\mathcal U}}
\def\SS{{\mathcal S}}
\def\DD{D}
\def\VV{{\mathcal V}}
\def\WW{{\mathcal W}}
\def\OO{{\mathcal O}}
\def\RR{{\mathcal R}}
\def\TT{{\mathcal T}}
\def\AA{{\mathcal A}}
\def\CC{{\mathcal C}}
\def\JJ{{\mathcal J}}
\def\NN{{\mathcal N}}
\def\HH{{\mathcal H}}
\def\XX{{\mathcal X}}
\def\XXX{{\mathscr X}}
\def\YY{{\mathcal Y}}
\def\ZZ{{\mathcal Z}}
\def\CC{{\mathcal C}}
\def\cir{{\scriptscriptstyle \circ}}
\def\circa{\thickapprox}
\def\vain{\rightarrow}
\def\parn{\par \noindent}
\def\salto{\vskip 0.2truecm \noindent}
\def\spazio{\vskip 0.5truecm \noindent}
\def\vs1{\vskip 1cm \noindent}
\def\fine{\hfill $\square$ \vskip 0.2cm \noindent}
\def\ffine{\hfill $\lozenge$ \vskip 0.2cm \noindent}
\newcommand{\rref}[1]{(\ref{#1})}
\def\beq{\begin{equation}}
\def\feq{\end{equation}}
\def\beqq{\begin{eqnarray}}
\def\feqq{\end{eqnarray}}
\def\barray{\begin{array}}
\def\farray{\end{array}}
%%%%%%%%% THIS NUMBERS EQUATIONS BY SECTIONS %%%%%%%%%%%%%
\makeatletter \@addtoreset{equation}{section}
\renewcommand{\theequation}{\thesection.\arabic{equation}}
%\thesection instead of \arabic{section} for correct equation numbering
% in appendices
\makeatother
%%%%%%%%%%%%%%%%%%%%%%%%%%%INTESTAZIONE%%%%%%%%%%%%%%%%%%%%%%%%%%%%%%%
\begin{titlepage}
{~}
\vspace{-2cm}
\begin{center}
{\huge On the constants in a Kato inequality for the Euler and Navier-Stokes equations}
\end{center}
\vspace{0.5truecm}
\begin{center}
{\large
Carlo Morosi$\,{}^a$, Livio Pizzocchero$\,{}^b$({\footnote{Corresponding author}})} \\
\vspace{0.5truecm} ${}^a$ Dipartimento di Matematica, Politecnico di Milano,
\\ P.za L. da Vinci 32, I-20133 Milano, Italy \\
e--mail: carlo.morosi@polimi.it \\
${}^b$ Dipartimento di Matematica, Universit\`a di Milano\\
Via C. Saldini 50, I-20133 Milano, Italy\\
and Istituto Nazionale di Fisica Nucleare, Sezione di Milano, Italy \\
e--mail: livio.pizzocchero@unimi.it
\end{center}
\begin{abstract}
We continue an analysis, started in
\cite{cok}, of some issues related
to the incompressible Euler or Navier-Stokes (NS) equations
on a $d$-di\-men\-sio\-nal torus $\Td$. More specifically,
we consider the quadratic term in these equations; this arises
from the bilinear map $(v, w) \mapsto v \sc \partial w$, where $v, w : \Td \vain \reali^d$ are
two velocity fields. We derive
upper and lower bounds for the constants in some
inequalities related to the above bilinear map; these bounds
hold, in particular, for the sharp constants
$G_{n d} \equiv G_n$ in the Kato inequality $| \la v \sc \partial w | w \ra_n |
\leqs G_n \| v \|_n \| w \|^2_{n}$, where
$n \in (d/2 + 1, + \infty)$ and
$v, w$ are in the Sobolev spaces $\HM{n}, \HM{n+1}$
of zero mean, divergence free vector fields of orders
$n$ and $n+1$, respectively.
As examples, the numerical values of our upper
and lower bounds are reported for $d=3$ and some values of $n$.
When combined with the results of \cite{cok} on another
inequality, the results of the present paper can be employed
to set up fully quantitative error estimates for
the approximate solutions of the Euler/NS equations, or
to derive quantitative bounds on the time of existence
of the exact solutions with specified initial data; a sketch of this program
is given.
\end{abstract}
\vspace{0.2cm} \noindent
\hfill \parn
\par \vspace{0.05truecm} \noindent \textbf{AMS 2000 Subject classifications:} 76D05, 26D10, 46E35.
\end{titlepage}
\section{Introduction}
The present paper continues our previous work on some inequalities
related to the Euler or Navier-Stokes (NS) equations. We work on a $d$-dimensional torus
$\Td$, and write these equations as
\beq {\partial u \over \partial t}  = - \LP(u \sc \partial u) + \nu \Delta u + f~, \label{eul} \feq
where: $u= u(x, t)$ is the divergence free velocity field; $x = (x_s)_{s=1,...,d} \in \Td$ are
the space variables (yielding the derivatives $\partial_s := \partial/\partial x_s$);
$\Delta := \sum_{s=1}^d \partial_{s s}$ is the Laplacian;
$(u \sc \partial u)_r := \sum_{s=1}^d u_s \partial_s u_r$ ($r=1,...,d$);
$\LP$ is the Leray projection onto the space of divergence free vector fields;
$\nu = 0$ for the Euler equations; $\nu \in (0,+\infty)$ (in fact $\nu=1$,
after rescaling) for the NS equations; $f = f(x,t)$
is the Leray projected density of external forces. As already noted
\cite{due}, the analysis of the above equations can be reduced to the case
where
the (spatial) means $\la u \ra := (2 \pi)^{-d} \int_{\Td} u \, d x$
and $\la f \ra$ are zero at all times. \parn
A precise functional setting for the above framework can be built using,
for suitable (integer or noninteger)
values of $n$, the Sobolev spaces
\beq \Hz{n}(\Td) \equiv \Hz{n} := \{ v : \Td \vain \reali^d~|~~
\sqrt{-\Delta}^{\,n} v \in \bb{L}^2(\Td),~ \la v \ra = 0 \}~, \feq
\beq \HM{n}(\Td) \equiv \HM{n} := \{ v \in \Hz{n}~|~\dive v=0 \} \feq
(the subscripts $0$, $\Sigma$ recall the
vanishing of the mean and of the divergence, respectively). For each $n$, we equip
$\Hz{n}$ with the standard inner product and the norm
\beq \la v | w \ra_n := \la \sqrt{-\Delta}^{\,n} v |  \sqrt{-\Delta}^{\,n} w \ra_{L^2}~,
\qquad \| v \|_n := \sqrt{\la v | v \ra_n}~, \feq
which can be restricted to the (closed) subspace $\HM{n}$. \parn
Our aim is to analyze quantitatively, in terms of the Sobolev inner products,
the quadratic map appearing in \rref{eul}. Some aspects of this map have been already
examined in the companion paper \cite{cok}; here we have considered the bilinear maps
sending two vector fields $v, w$ on $\Td$ into $v \sc \partial w$ or $\LP(v \sc \partial w)$,
and we have discussed some inequalities about them, the basic one being
\beq \| \LP(v \sc \partial  w) \|_n \leqs K_{n} \| v \|_{n} \| w \|_{n+1}
\qquad \mbox{for $n \in (\dd{d \over 2}, + \infty)$, $v \in \HM{n}$, $w \in \HM{n+1}$}~. \label{basineq} \feq
Our attention has been focused on the sharp constants
$K_{n} \equiv K_{n d}$ appearing therein, for which we have given
fully quantitative upper and lower bounds. \parn
In the present work we discuss other inequalities related to the quadratic Euler/NS nonlinearity,
discovered by Kato in \cite{Kato},
and establish upper and lower bounds for the unknown sharp
constants appearing therein. First of all we consider the inequality
\beq |\la v \sc \partial w | w \ra_n | \leqs G'_n \| v \|_n
\| w \|^{2}_{n} \quad \mbox{for $n \in (\dd{d \over 2} + 1, + \infty)$, $v \in \HM{n}$, $w \in \Hz{n+1}$}~,
\label{katinequa} \feq
writing $G'_{n} \equiv G'_{n d}$ for the sharp constants therein. With the additional assumption that
$w$ be divergence free, we can write
\beq |\la v \sc \partial w | w \ra_n | \leqs G_n \| v \|_n
\| w \|^{2}_{n} \quad \mbox{for $n \in (\dd{d \over 2} + 1, + \infty)$, $v \in \HM{n}$, $w \in \HM{n+1}$}~,
\label{katineq} \feq
with the sharp constant $G_{n} \equiv G_{n d}$ fulfilling the obvious relation $G_n \leqs G'_n$.
Let us observe that \rref{katineq} can be rephrased in terms of the Leray projection $\LP$; indeed,
with the assumptions therein we have $w = \LP w$ and this fact, combined with the symmetry
of $\LP$ in the Sobolev inner product, gives
\beq \la v \sc \partial w | w \ra_n = \la v \sc \partial w | \LP w \ra_n =
\la \LP(v \sc \partial w) | w \ra_n  \qquad
\mbox{for $v \in \HM{n}$, $w \in \HM{n+1}$}~. \label{18} \feq
Due to \rref{18}, Eq. \rref{katineq} is more directly related to the incompressible
Euler/NS equations \rref{eul}; in the sequel, \rref{katineq} is referred to as
the Kato inequality, and we call \rref{katinequa} the auxiliary Kato inequality. \parn
These inequalities (and similar ones)
are well known, but little has been done previously to evaluate with some accuracy the constants
which appear therein. On the other hand, quantitative bounds on such constants
are useful to estimate the time of existence of the solution of
\rref{eul} for a given initial datum, or its distance
from any approximate solution. \parn
In the present paper we derive
fully computable upper and lower bounds $G^{\pm}_n \equiv G^{\pm}_{n d}$ such that
\beq G^{-}_n \leqs G_n \leqs G'_{n} \leqs G^{+}_n \label{uplow} \feq
for all $n > d/2 + 1$. As examples, the bounds $G^{\pm}_n$ are computed in
dimension $d=3$, for some values of $n$. In these cases the upper and lower bounds
are not too far, at least for the purpose to apply them to the Euler/NS equations. \parn
To be more precise about such applications, let us exemplify
a framework already mentioned in \cite{cok}; the starting point of this setting is
a result of Chernyshenko, Constantin, Robinson and Titi \cite{Che}, that can be stated as
follows. Consider the Euler/NS equation \rref{eul} with a specified
initial condition $u(x,0) = u_0(x)$; let $u_{\ap} : \Td \times [0,T_\ap] \vain \reali^d$ be
an approximate solution of this Cauchy problem with errors $\er: \Td \times [0,T_\ap] \vain \reali^d$
on the equation and $\erd : \Td \vain \reali$ on the initial condition, by which we mean that
\beq \er := {\partial u_{\ap} \over \partial t} +  \LP(u_{\ap} \sc \partial u_{\ap}) - \nu \Delta u_{\ap} - f~,
\qquad
\erd := u_{\ap}(\cdot, 0) - u_0~. \feq
Fix $n \in (d/2+1, + \infty)$; then,
Eq. \rref{eul} with datum $u_0$ has a (strong) exact solution $u$ in $\HM{n}$
on a time interval $[0,T] \subset [0,T_\ap]$,
if $T$ and $u_{\ap}$ fulfill the inequality
\beq \| \erd \|_n + \int_{0}^T \| \er(t) \|_n d t < {1 \over G_n T}
\, e^{\dd{- \int_{0}^T (G_n \| u_{\ap}(t) \|_n +  K_n \| u_{\ap}(t) \|_{n+1}) d t }}
\label{criter} \feq
($u_{\ap}(t) := u_{\ap}(\cdot, t)$, $\er(t) := \er(\cdot, t)$).
For a given datum $u_0$, one can try a practical implementation of the above criterion
after choosing a suitable $u_{\ap}$ (say, a Galerkin approximate
solution). Of course, $T$ can be evaluated via \rref{criter} only in the presence of quantitative
information on $K_n$ and $G_n$, which are missing in \cite{Che}.
In a forthcoming paper \cite{coche}, our estimates on $K_n$ and $G_n$ will be
employed together with the existence condition \rref{criter} (or with some refinement of it,
suited as well to get bounds on $\| u(t) - u_{\ap}(t) \|_n$).
 \parn
For completeness we wish to mention that a program similar to the one described above,
but based on technically different inequalities, has been developed in \cite{due} \cite{tre}
for the incompressible NS equations in Sobolev spaces of lower order. For example,
in \cite{tre} we have considered the NS equations in $\HM{1}(\Tt)$; here we have derived
a fully quantitative upper bound on the vorticity $\| \curl u_0 \|_{L^2}$ of the initial
datum, which ensures global existence of the solution. \parn
Again for completeness, we remark that the fully quantitative attitude proposed
here for the Euler/NS equations is more or less close to the viewpoints
of other authors about these equations, or about different
nonlinear evolutionary PDEs \cite{Abd} \cite{Bart} \cite{Kyr}  \cite{Pel}
\cite{Zg1} \cite{Zg2}.
\salto
\textbf{Organization of the paper.}
Section \ref{inns} summarizes our standards about Sobolev spaces on $\Td$ and the Euler/NS
quadratic nonlinearity.
\parn
Section \ref{sequag} states
the main results of the paper; here we present our upper and lower bounds $G^{\pm}_n$ on
the constants in the inequalities \rref{katinequa} \rref{katineq},
which are treated by Propositions \ref{propupg} and \ref{prolowg}. The upper bounds
are determined by the $\sup$ of a positive function $\GG_n$,
defined on the space $\Zd \setminus \{0\}$ of nonzero Fourier wave vectors;
at each point $k \in \Zd \setminus \{0\}$,
$\GG_n(k)$ is a sum (of convolutional type) over
$\Zd \setminus \{0, k\}$. The lower bounds are determined by suitable trial functions.
As examples, in Eq. \rref{boug} we report the numerical values of
$G^{\pm}_{n}$, for $d=3$ and $n=3,4,5,10$. \parn
Section \ref{provewelg} contains the proofs of the previously mentioned
Propositions \ref{propupg}, \ref{prolowg}. \parn
Several appendices are devoted
to the practical evaluation of the function $\GG_n$ mentioned before,
and of the bounds $G^{\pm}_n$.
Appendix \ref{proveg}
presents some preliminary notations and results. Appendix \ref{appeg} contains
the main theorem (Proposition \ref{proggnd}) about the evaluation of $\GG_n$ and of its sup.
Appendix \ref{appe345g} gives details on the computation of $\GG_n$, and
on the corresponding upper bounds $G^{+}_n$,
for the previously mentioned cases $d=3$, $n=3,4,5,10$. Appendix \ref{appe345lowg}
describes the computation of the bounds $G^{-}_n$, for the same
values of $d$ and $n$.
\parn
For all the numerical computations required by this paper, as well as
for some lengthy symbolic manipulations,  we have used systematically
the software MATHEMATICA.
Throughout the paper,
an expression like $r= a. b c d e...~$ means the following: computation of
the real number $r$ via MATHEMATICA produces as an output $a.b c d e$,
followed by other digits not reported for brevity.
\section{Some preliminaries}
\label{inns}
We use for Sobolev spaces and the Euler/NS bilinear map
the same notations proposed in \cite{cok}; for the reader's convenience,
these are summarized hereafter.
Throughout the paper, we work in any space dimension
\beq d \geqs 2~; \feq
we use $r,s$ as indices running from $1$ to $d$.
For $a, b \in \complessi^d$ we put
\beq a \, \sc \, b := \sum_{r=1}^d a_r \, b_r~; \qquad |a| := \sqrt{\overline{a} \, \sc \, a} \feq
where $\overline{a} := (\overline{a_r})$ is the complex conjugate of $a$. We often refer
to the $d$-dimensional torus
\beq \Td := \underbrace{\To \times ... \times \To}_{\tiny{\mbox{$d$ times}}}~,
\qquad \To := \reali/(2 \pi \interi)~, \feq
whose elements are typically written
$x = (x_r)_{r=1,...d}$.
\salto
\textbf{Distributions on $\boma{\Td}$, Fourier series and
Sobolev spaces.}
The space of periodic distributions $\DD'(\Td, \complessi) \equiv \DD'_{\comple}$ is
the (topological) dual of $C^{\infty}(\Td, \complessi) \equiv C^{\infty}_{\comple}$;
$\la v, f \ra \in \complessi$ denotes the action of a distribution $v \in \DD'_{\comple}$
on a test function $f \in C^{\infty}_{\comple}$. \parn
Each $v \in \DD'_{\comple}$
has a unique (weakly convergent) Fourier series expansion
\beq v = \sum_{k \in \Zd} v_k e_k~, \quad
e_k(x) := {1 \over (2 \pi)^{d/2}} \, e^{i k \sc \, x}~\mbox{for $x \in \Td$}~, \quad
v_k := \la v, e_{-k} \ra~ \in \complessi~. \label{fs} \feq
The complex conjugate of a distribution $v \in \DD'_{\comple}$ is the unique distribution $\overline{v}$ such that
$\overline{\la v, f \ra} = \la \overline{v}, \overline{f} \ra$ for each $f
\in C^{\infty}_{\comple}$; one has $\overline{v}= \sum_{k \in \Zd} \overline{v_{k}} \, e_{-k}$. \parn
The \textsl{mean} of $v \in \DD'_{\comple}$ and the space
of \textsl{zero mean} distributions are
\beq \la v \ra := {1 \over (2 \pi)^d} \la v, 1 \ra = {1 \over (2 \pi)^{d/2}} v_0~,\qquad
\DD'_{\comple \zz} := \{ v \in \DD'_{\comple}~|~\la v \ra = 0 \}
\label{mean} \feq
(of course, $\la v, 1 \ra = \int_{\Td} v \, d x$ if $v \in L^1(\Td, \complessi, d x)$).
The relevant Fourier coefficients of zero mean distributions are labeled by the set
\beq \Zd_\zz := \Zd \setminus \{0\}~. \feq
The distributional derivatives
$\partial/\partial x_s \equiv \partial_s$ and the Laplacian
$\Delta := \sum_{s=1}^d \partial_{s s}$ send $\DD'_{\comple}$ into $\DD'_{\comple \zz}$ and,
for each $v$,
$\partial_s v = i \sum_{k \in \Zd_\zz} k_s v_k e_k$,
$\Delta v = - \sum_{k \in \Zd_\zz} | k |^2 v_k e_k$.
For any $n \in \reali$, we further define
\beq \sqrt{-\Delta}^{\, n} : \DD'_{\comple} \vain \DD'_{\comple \zz}~, \qquad
v \mapsto \sqrt{- \Delta}^{\, n}  v := \sum_{k \in \Zd_\zz} | k |^{n} v_k e_k~. \label{reg} \feq
The space of
\textsl{real} distributions is
\beq \DD'(\Td, \reali) \equiv \DD' := \{ v \in \DD'_{\comple}~|~\overline{v} = v \} =
\{ v \in \DD'_{\comple}~|~\overline{v_k} = v_{-k}~\mbox{for all $k \in \Zd$} \}~. \feq
For $p \in [1,+\infty]$ we often consider the real space
\beq L^p(\Td, \reali, d x) \equiv L^p~, \qquad \feq
\parn
especially for $p=2$. $L^2$ is a Hilbert space
with the inner product
$\la v | w \ra_{L^2} := \int_{\Td} v(x) w(x) d x = \sum_{k \in \Zd} \overline{v_k} w_k$ and
the induced norm
$\| ~ \|_{L^2}$. \parn
The zero mean parts of $\DD'$ and $L^p$ are
\beq \DD'_\zz := \{ v \in \DD'~|~~\la v \ra = 0 \}~,
\qquad  L^p_0 := L^p \cap \Dz~; \feq
all the differential operators mentioned before send $\DD'$ into $\DD'_{\zz}$. \parn
For each $n \in \reali$, the \textsl{zero mean Sobolev space}
$H^n_\zz(\Td, \reali) \equiv H^n_\zz$ is defined by
\beq  H^n_\zz := \{ v \in \DD'_\zz~|~\sqrt{-\Delta}^{\, n} v \in L^2 \}
= \{ v \in \DD'_\zz~|~\sum_{k \in \Zd_\zz} | k |^{2 n} | v_k |^2 < + \infty~\}~;
\label{hn} \feq
this is a real Hilbert space with the inner product
$\la v | w \ra_{n} := \la \sqrt{-\Delta}^{\, n} v \, | \, \sqrt{- \Delta}^{\, n} w \ra_{L^2}$
$= \sum_{k \in \Zd_\zz} | k |^{2 n} \, \overline{v_k}  w_k$ and the induced norm
$\|~ \|_{n}$. Of course, $H^0_0 = L^2_0$.
\salto
\textbf{Spaces of vector valued functions on $\Td$.}
If $\Vi(\Td, \reali) \equiv \Vi$ is any vector space of \textsl{real}
functions or distributions on $\Td$, we write
\beq \bVi(\Td) \equiv \bVi := \{ v = (v_1,...,v_d)~|~v_r \in \Vi \quad \mbox{for all $r$}\}~. \feq
In this way we can define, e.g., the spaces $\bb{\DD}'(\Td) \equiv \bb{\DD}'$,
$\bb{L}^p(\Td) \equiv \bb{L}^p$ ($p \in [1,+\infty]$), $\Hz{n}(\Td) \equiv \Hz{n}$.
Any $v = (v_r) \in \bb{\DD}'$ is referred to
as a (distributional) \textsl{vector field} on $\Td$. We note that $v$ has a unique Fourier
series expansion \rref{fs} with coefficients
\beq v_k := (v_{r k})_{r =1,...,d} \in \complessi^d~, \qquad v_{r k}  := \la v_r, e_{-k} \ra~; \feq
as in the scalar case, the reality of $v$ ensures $\overline{v_k} = v_{-k}$. \parn
$\bb{L}^2$ is a real Hilbert space, with the inner product and the norm
\beq \la v | w \ra_{L^2} := \int_{\Td} v(x) \sc w(x) d x = \sum_{k \in \Zd} \overline{v_k} \sc w_k~,
\qquad \| v \|_{L^2} := \sqrt{\la v | v \ra_{L^2}}~. \feq
We define componentwise the mean $\la v \ra \in \reali^d$
of any $v \in \bb{\DD}'$ (see Eq. \rref{mean});
$\Dz$ is the space of zero mean vector fields, and
$\bb{L}^p_\zz = \bb{L}^p \cap \Dz$. \parn
We similarly define
componentwise the operators
$\partial_s, \Delta, \sqrt{-\Delta}^{\,n} : \bb{\DD}' \vain \Dz$. \parn
For any real $n$, the $n$-th Sobolev space of zero mean vector fields $\Hz{n}(\Td) \equiv \Hz{n}$ is made
of all $d$-uples $v$ with components
$v_r \in H^n_\zz$; an equivalent definition can be given via Eq.\rref{hn},
replacing therein $L^2$ with $\bb{L}^2$.
$\Hz{n}$ is a real Hilbert space with the inner product and the induced norm
\parn
\vbox{
\beq
\la v | w \ra_{n} := \la \sqrt{-\Delta}^{\, n} \, v \, | \sqrt{-\Delta}^{\,n} \,w \ra_{L^2}
= \sum_{k \in \Zd_\zz} | k |^{2 n} \, \overline{v_k}  \sc \, w_k~, \feq
$$ \| v \|_{n} = \| \sqrt{-\Delta}^{\,n} v \|_{L^2} =
\sqrt{ \sum_{k \in \Zd_\zz} | k |^{2 n} \, |v_k|^2 }~. $$
}
\noindent
\salto
\textbf{Divergence free vector fields.}
Let
$\dive : \bb{\DD}' \vain \DD'_\zz$, $v \mapsto \dive v := \sum_{r=1}^d \partial_r v_r$
$= i \, \sum_{k \in \Zd_\zz}  (k \sc \, v_k) e_k$. Hereafter we introduce
the space $\Ds$ of \textsl{divergence free (or solenoidal) vector fields} and some
subspaces of it, putting
\beq \Ds := \{ v \in \bb{\DD}'~|~\dive v = 0 \}
= \{ v \in \bb{\DD}'~|~k \sc \, v_k = 0~\forall k \in \Zd~\} ~; \feq
\beq ~~\Dsz := \Ds \cap \Dz~, \quad \bb{L}^p_{\ss} := \bb{L}^p \cap \Ds~,~\bb{L}^p_{\so} := \bb{L}^p
\cap \Dsz \quad (p \in [1,+\infty])~, \feq
\beq \HM{n} := \Ds \cap \Hz{n} \quad (n \in \reali). \feq
$\HM{n}$ is a closed subspace of the Hilbert space $\Hz{n}$,
that we equip with the restrictions of  $\la~|~\ra_n$, $\|~\|_n$.
The \textsl{Leray projection} is the (surjective) map
\beq \LP : \bb{\DD}' \vain \Ds~, \qquad
v \mapsto \LP v := \sum_{k \in \Zd} (\LP_k v_k) e_k~, \label{ler1} \feq
where, for each $k$,
$\LP_k$ is the orthogonal projection of $\complessi^d$ onto the orthogonal
complement of $k$;
more explicitly, if $c \in \complessi^d$,
\beq \LP_\zz c = c~, \qquad \LP_k c = c - {k \sc \, c  \over | k |^2}\, k \quad \mbox{for $k \in \Zd_\zz$}~.
\label{ler2} \feq
From the Fourier representations of $\LP$, $\la~\ra$, etc., one easily infers
\beq \la \LP v \ra = \la v \ra~\mbox{for $v \in \bb{\DD}'$},
\quad \LP \Dz = \bb{D}'_{\so}, \quad \LP \bb{L}^2= \bb{L}^2_{\ss},
\quad \LP \Hz{n} = \HM{n}~\mbox{for}~ n \in \reali~. \label{tvv} \feq
Furthermore, $\LP$ is an orthogonal projection in each one of the Hilbert spaces $\bb{L}^2$, $\Hz{n}$; in
particular,
\beq \| \LP v \|_n \leqs \| v \|_n \qquad \mbox{for $v \in \Hz{n}$}~. \feq
\salto
\textbf{Making contact with the Euler/NS equations.}
The quadratic nonlinearity in the Euler/NS equations is related to
the bilinear map sending two (sufficiently regular) vector fields
$v, w$ on $\Td$ into $v \sc \partial w$; we are now ready to discuss this map. \parn
Hereafter we often refer to the case
\beq v \in \bb{L}^2~, \quad \partial_s w \in \bb{L}^2~~ (s =1,...,d)~; \label{asin} \feq
the above condition on the derivatives of $w$ implies $w \in \bb{L}^2$. \parn
The results mentioned in the sequel are known: the proofs of Lemmas
\ref{lemma1}, \ref{lemma3} are found, e.g., in \cite{cok},
and the proof of Lemma \ref{propzero} is reported only for completeness.
\begin{prop}
\label{lemma1}
\textbf{Lemma.} For $v, w$ as in \rref{asin}, consider the vector field
$v \sc \partial w$ on $\Td$, of components
\beq (v \sc \partial w)_r : = \sum_{s=1}^d v_s \partial_s w_r~; \feq
this is well defined and belongs to $\bb{L}^1$.
With the additional assumption $\dive v = 0$,
one has $\la v \sc \partial w \ra = 0$ (which also implies $\la \LP(v \sc \partial w) \ra = 0$, see
\rref{tvv}).
\end{prop}
\begin{prop}
\label{lemma3}
\textbf{Lemma.} Assuming \rref{asin}, $v \sc \partial w$ has Fourier coefficients
\beq (v \sc \partial w)_k = {i \over (2 \pi)^{d/2}} \sum_{h \in \Zd} [v_{h} \sc \, (k - h)] w_{k - h}
\qquad \mbox{for all $k \in \Zd$}~. \label{infert} \feq
\end{prop}
\begin{prop}
\label{propzero}
\textbf{Lemma.} Besides \rref{asin}, assume
$\dive  v=0$ and $v \sc \partial w \in \bb{L}^2$. Then
\beq \la v \sc \partial w |w \ra_{L^2} = 0~. \label{tethen} \feq
\end{prop}
\textbf{Proof.} Suppose for a moment that $v,w : \Td \vain \reali^d$ are $C^1$, with no other condition; then
(integrating by parts in one passage)
\parn
\vbox{
$$\la v \sc \partial w | w \ra_{L^2}= \sum_{r,s=1}^d \int_{\Td} v_s (\partial_s w_r)  w_r \, d x
={1 \over 2} \sum_{r,s=1}^d \int_{\Td} v_s \partial_s(w_r^2)\, d x $$
$$ = - {1 \over 2} \sum_{r,s=1}^d \int_{\Td} (\partial_s v_s) w_r^2 \, d x
= - {1 \over 2} \int_{\Td} (\dive v) |w |^2 \, d x~. $$
}
\noindent
In particular, \rref{tethen} holds if $v,w$ are $C^1$ and $\dive v = 0$.
By a density argument, one extends \rref{tethen} to all $v,w$ as in the statement of the
Lemma. \fine
The following result, essential for the sequel, is also well known (see, e.g., \cite{cok}).
\begin{prop}
\label{start}
\textbf{Proposition.} Let $n \in (d/2, + \infty)$. If
$v \in \HM{n}$ and $w \in \Hz{n+1}$, one has $v \sc \partial w \in \Hz{n}$.
Furthermore, the map $(v,w) \mapsto v \sc \partial w$
is bilinear and continuous between the spaces mentioned before.
\end{prop}
\section{The Kato inequality}
\label{sequag}
Throughout this section we assume
\beq n \in ({d \over 2} + 1, + \infty)~. \label{assdg} \feq
The following Proposition \ref{weldefg} is known, dating back to \cite{Kato} (see \cite{CoFo} for
a more general formulation, similar to the one proposed hereafter). As a matter of fact,
the quantitative analysis presented later in this paper
also gives, as a byproduct, an alternative proof of this Proposition.
\begin{prop}
\label{weldefg}
\textbf{Proposition.}
Let $v \in \HM{n}$, $w \in \Hz{n+1}$ (so that $v \sc \partial w \in \Hz{n}$).
Then, there is $G' \in [0,+\infty)$, independent of $v,w$, such that
\beq | \la v \sc \partial w |w  \ra_n | \leqs G' \| v \|_n \| w \|^2_n~. \label{diseqg} \feq
\end{prop}
\begin{prop}
\textbf{Definition.} We put
\beq G'_{n d} \equiv G'_n  \label{gpnd} \feq
$$ := \min \{ G' \in [0,+\infty)~|~|\la v \sc \partial w  | w \ra_n | \leqs G' \| v \|_n
\| w \|^{2}_{n}~\mbox{for all $v \in \HM{n}$, $w \in \Hz{n+1}$}\}~; $$
\beq G_{n d} \equiv G_n  \label{gnd} \feq
$$ := \min \{ G \in [0,+\infty)~|~|\la v \sc \partial w  | w \ra_n | \leqs G \| v \|_n
\| w \|^{2}_{n}~\mbox{for all $v \in \HM{n}$, $w \in \HM{n+1}$}\}~. $$
(Note that all $w$'s in \rref{gnd} are divergence free, a property not required in
\rref{gpnd}.)
\end{prop}
With the language of the Introduction, $G'_n$ and $G_n$ are, respectively,
the sharp constants in the ''auxiliary Kato inequality'' \rref{katinequa}
and in the Kato inequality \rref{katineq}; we recall that $\LP$ could be
inserted into \rref{gnd}, due to the relation \rref{18}
$\la v \sc \partial w  | w \ra_n = \la \LP(v \sc \partial w)  | w \ra_n$.
It is obvious that
\beq G_n \leqs G'_n \,; \label{ensure} \feq
in the rest of the section (which is its original part) we present
computable upper and lower bounds on $G'_n$ and $G_{n}$, respectively. \parn
The upper bound requires a more lengthy analysis; the final result relies on a
function $\GG_{n d} \equiv \GG_{n}$, appearing in the forthcoming Definition \ref{deggnd}.
To build this function, as in \cite{cok} we refer to the exterior power $\We^2 \reali^d$,
identified with the space of real, skew-symmetric $d \times d$ matrices
$A = (A_{r s})_{r,s,=1,...,d}$. We
consider the (bilinear, skew-symmetric) operation $\we$ and the norm $|~|$ defined by
\beq \we : \reali^d \times \reali^d \vain {\We}^2 \reali^d,
\quad (p,  q) \mapsto p \we q~~\mbox{s.t.}~~
(p \we q)_{r s} := p_r q_s - q_r p_s~;
\label{ester} \feq
\beq |~| : {\We}^2 \reali^d \vain [0,+\infty),
\qquad A = (A_{r s}) \mapsto |A| :=  \sqrt{{1 \over 2}\, \sum_{r,s=1}^d  |A_{r s}|^2}~.
\label{normest}
\feq
In the sequel, for $p,q \in \reali^d$, we often use the relations
\beq |p \we q| = \sqrt{|p|^2 |q|^2 - (p \sc q)^2} = |p| |q| \sin \tet~, \label{norqp} \feq
where $\tet \equiv \tet(p, q) \in [0,\pi]$ is the convex angle between $p$
and $q$ (defined arbitrarily, if $p=0$ or $q=0$); we use as well the inequality
\beq |p \we q|  \leqs |p||q|~. \feq
Keeping in mind these facts, let us stipulate the following.
\begin{prop}
\label{deggnd}
\textbf{Definition.} We put
\beq \Zd_{0 k} := \Zd \setminus \{0,k\} \qquad \mbox{for each $k \in \Zd_{0}$}~; \feq
\beq \GG_{n d } \equiv \GG_n : \Zd_0 \vain (0,+\infty),\quad
k \mapsto \GG_{n}(k) :=
\sum_{h \in \Zd_{0 k}} {|h \we k|^2 (|k|^n - |k-h|^n)^2
\over |h|^{2 n + 2} |k-h|^{2 n}}~. \label{ggnd}\feq
\end{prop}
\begin{rema}
\label{rem34}
\textbf{Remarks.}
(i) For any $k \in \Zd_0$ one has $\GG_{n}(k) < + \infty$, as stated above, since
\beq  {|h \we k|^2 (|k|^n - |k-h|^n)^2 \over |h|^{2 n + 2} |k-h|^{2 n}} = O({1 \over |h|^{2 n}}) \qquad
\mbox{for $h \vain \infty$}~, \feq
and $2 n > d$. \parn
(ii) Consider the reflection operators
$R_r(k_1,.., k_r,...,k_d)$ $:= (k_1,...,$ $-k_r,...,k_d)$
($r=1,...,d$)
and the permutation operators $P_{\si}(k_1,...,k_d) :=
(k_{\si(1)},...,k_{\si(d)})$ ($\si$ a permutation of $\{1,...,d\}$); then
\beq \GG_{n}(R_r k) = \GG_{n}(k)~, \quad \GG_{n}(P_\sigma k) = \GG_{n}(k)~
\qquad \mbox{for each $k \in \Zd_0$}~. \label{clai} \feq
The proof
is very similar to the one employed for the analogous properties
of the function $\KK_{n}$ appearing in \cite{cok}.
\parn
(iii) In Appendix \ref{appeg} we will prove that
\beq \sup_{k \in \Zd_0} \GG_{n}(k) < + \infty~, \label{supfinito} \feq
and give tools for the practical evaluation of $\GG_{nd}$ and of its sup.
\fine
\end{rema}
The main result of the present section is the following.
\begin{prop}
\label{propupg}
\textbf{Proposition.}
The constant $G'_{n}$ defined by \rref{gnd} has the upper bound
\beq G'_n \leqs G^{+}_n~, \label{hasup} \feq
\beq G_{n} := {1 \over (2 \pi)^{d/2}}
\sqrt{\sup_{k \in \Zd_0} \GG_{n}(k)}~\mbox{(or any approximant for this)}. \label{desup} \feq
\end{prop}
\textbf{Proof.} See Section \ref{provewelg}. \fine
The practical calculation of the above upper bound is made possible
by a general method, illustrated in Appendix \ref{appeg};
the results of such calculations, for $d=3$ and some
illustrative choices of $n$,
are reported at the end of this section. \parn
Let us pass to the problem of finding a lower bound for
the constant $G_{n}$; this can be obtained
directly from the tautological inequality
\beq G_n \geqs {|\la v \sc \partial w | w \ra_n| \over \| v \|_n \| w \|^2_{n}}
\qquad \mbox{for $v \in \HM{n} \setminus \{0\}$, $w \in \HM{n+1} \setminus \{0\}$}~,
\label{gneq} \feq
choosing for $v$ and $w$ two suitable
non zero ``trial functions''; hereafter we consider a choice
where
$v_k = 0$ for $k \in \Zd_0 \setminus V$ and $w_k = 0$ for $k \in \Zd_0 \setminus W$
with $V,W$ two finite sets. For the sake of brevity in the exposition of the final result,
let us stipulate the following.
\begin{prop} \textbf{Definition.}
We put
\beq \S_d \equiv \S := \{ (u_k)_{k \in U}~|~U \subset \Zd_0~\mbox{finite}, -U=U; \label{esseu} \feq
$$ u_k \in \complessi^d,~ \overline{u_{k}} = u_{-k},~k \sc u_k = 0~\mbox{for all $k \in U$} \} $$
(the set $U$ can depend on the family $(u_k)$, and $-U := \{ - k~|~k \in U \}$). \fine
\end{prop}
\begin{prop}
\label{prolowg}
\textbf{Proposition.}
Consider two nonzero families $(v_k)_{k \in V}$, $(w_k)_{k \in W}$ $\in \S$;  these
give the lower bound
\beq G_n \geqs G^{-}_n~, \label{gnlowlow} \feq
where
\beq G^{-}_n := {1 \over (2 \pi)^{d/2}} \, {|P_n((v_k), (w_k))| \over N_n((v_k)) N^2_n((w_k))}~
\mbox{(or any lower approximant for this)}~,
\label{gnlow} \feq
$$ N_n((v_k)) := \Big(\sum_{k \in V} |k|^{2 n} |v_k|^2 \Big)^{1/2}, \qquad
N_n((w_k)) := \Big(\sum_{k \in V} |k|^{2 n} |w_k|^2 \Big)^{1/2}~, $$
$$ P_n((v_k), (w_k)) := - i \!\!\!\! \sum_{h \in V, \ell \in W, h + \ell \in W}
\!\!\!\! |h + \ell|^{2 n} (\overline{v_h} \sc \ell) (\overline{w_{\ell}} \sc w_{h + \ell})~. $$
\end{prop}
\textbf{Proof.} See Section \ref{provewelg}. Here, we anticipate the main idea:
the vector fields $v := \sum_{k \in V} v_k e_k$, $w := \sum_{k \in W} w_k e_k$ belong to
$\HM{m}$ for each real $m$, and $\| v \|_{n} = N_n((v_k))$,
 $\| w \|_{n} = N_n((w_k))$, $\la v \sc \partial w | w \ra_n =  (2 \pi)^{-d/2} P_n((v_k), (w_k))$;
so, \rref{gnlowlow} is just the relation \rref{gneq} for this choice of $v, w$.
\fine
Putting together Eqs. \rref{ensure} \rref{hasup} \rref{gnlowlow} we obtain
a chain of inequalities, anticipated in the Introduction,
$$ G^{-}_n \leqs G_n \leqs G'_{n} \leqs G^{+}_n~; $$
here, the bounds $G^{\pm}_n$ can be computed explicitly from their definitions
\rref{desup} \rref{gnlow}.
\begin{rema}
\textbf{Examples.}
For $d=3$ and $n=3,4,5,10$, Eq. \rref{desup} and Eq. \rref{gnlow}
(with suitable choices of $(v_k)$, $(w_k)$),
give
\beq G^{-}_3 = 0.114 \,,~ G^{+}_{3} = 0.438 \,; \qquad G^{-}_4 = 0.181 \,,~ G^{+}_{4} = 0.484 \,;
\label{boug} \feq
$$ G^{-}_5 = 0.280 \,,~ G^{+}_{5} = 0.749 \,; \qquad
G^{-}_{10} = 2.41 \,,~G^{+}_{10} = 7.56 $$
(see Appendices \ref{appe345g} and \ref{appe345lowg} for the upper and lower
bounds, respectively).
In the above, the ratios $G^{-}_{n}/G^{+}_n$
are $0.260...$, $0.373...$, $0.373...$, $0.318...$
for $n=3,4,5,10$, respectively. To avoid misunderstandings related to these examples,
we repeat that the approach of this paper
applies as well to noninteger values of $n$.
\end{rema}
\section{Proof of Propositions (\ref{weldefg} and) \ref{propupg}, \ref{prolowg}}
\label{provewelg}
For the reader's convenience, we report a Lemma from \cite{cok}.
\begin{prop}
\label{lemabz}
\textbf{Lemma.} Let
\beq p, q \in \reali^d \setminus \{0 \}~, \quad z \in \complessi^d,~p \sc z = 0~, \feq
and $\tet(p,q) \equiv \tet \in [0,\pi] $ be the convex angle between $q$ and $p$. Then
\beq | q \sc z | \leqs \sin \tet \,  |q | |z | = {|p \we q| \over |p|} \, |z |~.
\label{qscz} \feq
\end{prop}
\noindent
From now on,
$n \in (\dd{d \over 2} + 1,+\infty)$.
Hereafter we present an argument proving
(Proposition \ref{weldefg} and, simultaneously)
Proposition \ref{propupg}. This is divided in several steps; in particular,
Step 1 relies on an idea of Constantin and Foias \cite{CoFo}. These
authors use their idea to obtain a proof of the Kato inequalities, but
are not interested in the quantitative evaluation of the sharp constants
therein; our forthcoming argument can be regarded as a refined,
fully quantitative version of their approach, developed for
the specific purpose to estimate $G'_n$.
\salto
\textbf{Proof of Propositions \ref{weldefg}, \ref{propupg}}.
We choose $v \in \HM{n}$, $w \in \Hz{n+1}$ and proceed in some steps. \parn
\textsl{Step 1. We have
$v \in \bb{L}^{\infty}_{\so}$,
$\sqrt{-\Delta}^{\,n} w
\in \Hz{1}$, $v \sc \partial (\sqrt{-\Delta}^{\,n} w) \in \bb{L}^{2}_{\zz}$,
$v \sc \partial w \in \Hz{n}$ and $\sqrt{-\Delta}^{\,n} (v \sc \partial w) \in \bb{L}^{2}_{\zz}$;
furthermore, the vector field}
\beq z := \sqrt{-\Delta}^{\,n} (v \sc \partial w) - v \sc \partial (\sqrt{-\Delta}^{\,n} w) \in \bb{L}^{2}_{\zz}
\label{defze}\feq
\textsl{fulfills the equality}
\beq  \la v \sc \partial w |w \ra_{n} = \la z | \sqrt{-\Delta}^{\,n} w \ra_{L^2}~, \label{wehave} \feq
\textsl{which implies}
\beq | \la v \sc \partial w |w \ra_{n} | \leqs \| z \|_{L^2} \| w \|_n~. \label{wehave2} \feq
To prove all this, we first recall the Sobolev imbedding
$H^n_{\zz} \subset L^{\infty}$, holding because $n > d/2$ (see,
e.g., \cite{Ada}); this obviously implies
$\HM{n} \subset \bb{L}^{\infty}_{\so}$, so $v \in \bb{L}^{\infty}_{\so}$.
Of course, $\sqrt{-\Delta}^{\,n}$ sends $\Hz{n+1}$ into $\Hz{1}$,
thus $\sqrt{-\Delta}^{\,n} w \equiv u \in \Hz{1}$.
This implies $\partial_s u_{r} \in L^2$ that, with $v_s \in L^{\infty}$, gives
$(v \sc \partial u)_r= \sum_{s=1}^d v_s \partial_s u_{r} \in L^2$. Summing up,
$v \sc \partial u \in \bb{L}^2$; furthermore, $v \sc \partial u \in \bb{L}^2_0$
due to Lemma \ref{lemma1}.
The statement
$v \sc  \partial w \in \Hz{n}$ holds due to Proposition \ref{start}; since $\sqrt{-\Delta}^{\,n}$
sends $\Hz{n}$ into $\Hz{0} = \bb{L}^2_{\zz}$, we finally obtain
$\sqrt{-\Delta}^{\,n}(v \sc  \partial w)  \in \bb{L}^2_{\zz}$. \parn
To go on, we note that
\parn
\vbox{
$$ \la v \sc  \partial w|w \ra_{n} = \la \sqrt{-\Delta}^{\,n} ( v \sc  \partial w) | \sqrt{-\Delta}^{\,n} w \ra_{L^2} $$
$$ = \la \sqrt{-\Delta}^{\,n} (v \sc  \partial w) -
v \sc \partial( \sqrt{-\Delta}^{\,n} w) | \sqrt{-\Delta}^{\,n} w \ra_{L^2}
= \la z | \sqrt{-\Delta}^{\,n} w \ra_{L^2}~.
$$} \noindent
In the above: the first equality corresponds to the definition of $\la~|~\ra_n$,
the second one holds because
$\la v \sc \partial (\sqrt{-\Delta}^{\,n} w) | \sqrt{-\Delta}^{\,n} w \ra_{L^2} = 0$ by Lemma
\ref{propzero} (here applied to the vector fields $v,\sqrt{-\Delta}^{\,n} w$);
the last equality corresponds to the definition of $z$, and proves
Eq. \rref{wehave}. Now, the Schwartz inequality yields
$| \la v \sc \partial w |w \ra_{n} |$ $\leqs \| z \|_{L^2} \| \sqrt{-\Delta}^{\,n} w \|_{L^2}$
$= \| z \|_{L^2} \| w \|_n$, as in \rref{wehave2}.
\parn
\textsl{Step 2. The vector field $z$ in \rref{defze} has
Fourier coefficients
\beq z_{k} = {i \over (2 \pi)^{d/2}}
\sum_{h \in \Zd_{0 k}} [v_{h} \sc \, (k - h)] (|k|^ n - |k - h|^n) w_{k - h}
\quad \mbox{for all $k \in \Zd_0$}~.
\label{znk} \feq}
To prove this, let us start from the Fourier coefficients of
$v \sc \partial w$; this has zero mean,
so $(v \sc \partial w)_0=0$. The other coefficients are
\beq (v \sc \partial w)_k = {i \over (2 \pi)^{d/2}} \sum_{h \in \Zd_{0 k}} [v_{h} \sc \, (k - h)] w_{k - h}
\qquad \mbox{for all $k \in \Zd_0$}~; \label{infertt} \feq
this follows from \rref{infert} taking into account that, in the sum therein,
the term with $h=0$ vanishes due to $v_0=0$, and the term with $h =k$
is zero for evident reasons. \parn
Consider any $k \in \Zd_0$; Eq. \rref{infertt} implies
\beq [\sqrt{-\Delta}^{\,n} (v \sc \partial w) ]_k = |k|^n (v \sc \partial w)_k =
{i |k|^n \over (2 \pi)^{d/2}} \sum_{h \in \Zd_{0 k}} [v_{h} \sc \, (k - h)] w_{k - h}~. \label{dasub} \feq
The analogue of Eq. \rref{infertt} for the pair $v, \sqrt{-\Delta}^{\,n} w$ reads
$[v \sc \partial(\sqrt{-\Delta}^{\,n} w)]_k
= i (2 \pi)^{-d/2}$ $ ~\sum_{h \in \Zd_{0 k}} [v_{h} \sc \, (k - h)]
(\sqrt{-\Delta}^{\,n} w)_{k - h}$, i.e.,
\beq [v \sc \partial (\sqrt{-\Delta}^{\,n} w)]_k =
{i \over (2 \pi)^{d/2}} \sum_{h \in \Zd_{0 k}} [v_{h} \sc \, (k - h)] |k - h|^n w_{k - h}~.
\label{dasub2} \feq
Subtracting \rref{dasub2} from \rref{dasub}, we obtain the thesis \rref{znk}.
\parn
\textsl{Step 3. Estimating the Fourier coefficients of $z$}.
Let $k \in \Zd_0$; Eq. \rref{znk} implies
\beq |z_k| \leqs {1 \over (2 \pi)^{d/2}}
\sum_{h \in \Zd_{0 k}} |v_{h} \sc \, (k - h)|~\Big||k|^ n - |k - h|^n\Big|~ |w_{k - h}|~.
\label{zznk} \feq
To go on, we note that
$h \sc v_h = 0$ due to the assumption $\dive v=0$;
so, we can apply Eq. \rref{qscz} with $p=h$, $q=k-h$ and $z = v_h$, which gives
\beq |v_{h} \sc \, (k - h)| \leqs {|h \we (k-h)| \over |h|} |v_h| =
{|h \we k| \over |h|} |v_h| \label{ins1} \feq
(recall that $h \we (k-h) = h \we k$).
Inserting the inequality \rref{ins1} into \rref{zznk}, we get
\parn
\vbox{
\beq |z_k| \leqs {1 \over (2 \pi)^{d/2}}
\sum_{h \in \Zd_{0 k}} {|h \we k| \over |h|} |v_h| ~\Big|\,|k|^ n - |k - h|^n \Big| \, |w_{k - h}| \feq
$$ = {1 \over (2 \pi)^{d/2}} \sum_{h \in \Zd_{0 k}}
{|h \we k| \,\Big||k|^ n - |k - h|^n\,\Big| \over |h|^{n + 1} |k - h |^{n}}
\Big( |h |^n | v_h | |k - h|^{n} | w_{k - h} | \Big)~. $$
}
\noindent
Now, H\"older's inequality $| \sum_h~ a_h b_h |^2 \leqs \Big(\sum_h | a_h |^2\Big)
\Big(\sum_h~| b_h |^2 \Big)$ gives
\parn
\vbox{
\beq | z_k |^2 \leqs {1 \over (2 \pi)^d} \, \GG_{n}(k) \QQ_{n}(k)~\mbox{for
all $k \in \Zd_0$}, \label{dainsg} \feq
$$ \GG_{n}(k) :=  \sum_{h \in \Zd_{0 k}} {|h \we k|^2
(|k|^ n - |k - h|^n)^2 \over | h |^{2 n + 2} | k - h |^{2 n} } \quad \mbox{as in \rref{ggnd}}, $$
$$ \QQ_{n}(k) \equiv \QQ_{n}(v,w)(k) := \sum_{h \in \Zd_{0 k}} |h|^{2 n} |  v_{h} |^2
| k - h |^{2 n} | w_{k - h} |^2 $$
}
\noindent
(in the definition of $\QQ_{n}(k)$ one can write as well $\sum_{h \in \Zd_{0}}$, since
the general term of the sum vanishes for $h=k$).
\parn
\textsl{Step 4. Estimates on $\| z \|_{L^2}$}.
Eq. \rref{dainsg} implies
$$ \| z_{n} \|^2_{L^2} = \sum_{k \in \Zd_0} | z_k |^{2} \leqs
{1 \over (2 \pi)^d} \sum_{k \in \Zd_0} \GG_{n}(k) \QQ_{n}(k)~\leqs
{1 \over (2 \pi)^d} \Big(\sup_{k \in \Zd_0} \GG_{n}(k)\Big)
\Big(\sum_{k \in \Zd_0} \QQ_{n}(k)\Big)~. $$
The sup of $\GG_{n}$ is finite, as we will show (by an independent argument) in
Proposition \ref{proggnd};
making reference to the definition of $G^{+}_n$ in terms of this sup (see Eq. \rref{desup}),
we can write the last result as
\beq \| z_{n} \|^2_{L^2} \leqs (G^{+}_n)^2 \sum_{k \in \Zd_0} \QQ_{n}(k)~. \label{daretg} \feq
On the other hand,
$$ \sum_{k \in \Zd_0} \QQ_{n}(k) =
\sum_{h \in \Zd_0} |h|^{2 n} |  v_{h} |^2 \sum_{k \in \Zd_{0}} | k - h |^{2 n} | w_{k - h} |^2
= \sum_{h \in \Zd_0} |h|^{2 n} |  v_{h} |^2 \sum_{\ell \in \Zd_{0 h}} | \ell |^{2 n} | w_{\ell} |^2  $$
\beq  \leqs \sum_{h \in \Zd_0} |h|^{2 n} |  v_{h} |^2~
\sum_{\ell \in \Zd_{0}} | \ell |^{2 n} | w_{\ell} |^2  = \| v \|^2_{n} \| w \|^2_{n}~. \feq
Inserting this result into \rref{daretg}, we obtain
\beq \| z \|_{L^2} \leqs G^{+}_n
\| v \|_{n} \| w \|_{n}~. \label{tellsg} \feq
\textsl{Step 5. Concluding the proofs of Propositions \ref{weldefg}, \ref{propupg}.}
Eqs. \rref{wehave2} \rref{tellsg} imply
\beq | \la v \sc \partial w |w \ra_{n} | \leqs G^{+}_n \, \| v \|_{n} \| w \|^{2}_n~;
\label{imply} \feq
so, Proposition \ref{weldefg} is proved. Eq. \rref{imply} also indicates that
the sharp constant $G'_n$ in \rref{gpnd} fulfills $G'_{n} \leqs G^{+}_n$;
this proves Eq. \rref{hasup} and Proposition \ref{propupg}. \fine
\salto
We conclude this section proving the statements of Section \ref{sequag} on the lower bounds $G^{-}_n$. \salto
\textbf{Proof of Proposition \ref{prolowg}.} Let us recall the definition \rref{esseu} of
$\S$;
our argument is divided in some steps. \parn
\textsl{Step 1. Let $(u_k)_{k \in U} \in \S$. Then,
\beq u := \sum_{k \in U} u_k e_k \feq
belongs to $\HM{m}$ for each real $m$, and
\beq \| u \|_m = \left( \sum_{k \in U} |k|^{2 m} |u_k|^2 \right)^{1/2} \equiv N_m((u_k))~. \feq}
These statements are self-evident; of course, the conditions $\overline{u_{k}} = u_{-k}$ and
$k \sc u_k = 0$ in \rref{esseu} ensure $u$ to be real, and divergence free. \parn
\textsl{Step 2. Consider two families $(v_k)_{k \in V}$, $(w_k)_{k \in W}$ $\in \S$, and
define $v := \sum_{k \in V} v_k e_k$, $w := \sum_{k \in W} w_k e_k$ . Then
\beq \la v \sc \partial w | w \ra_n =  {1 \over (2 \pi)^{d/2}} P_n((v_k), (w_k)) \label{lappra} \feq
where, as in \rref{gnlow}, $P_n(v,w) :=  - i \sum_{h \in V, \ell \in W, h + \ell \in W}
|h + \ell|^{2 n} (\overline{v_h} \sc \ell) (\overline{w_{\ell}} \sc w_{h + \ell})$.} \parn
In fact, the Fourier coefficients of $v \sc \partial w$ have the expression \rref{infert}
$$ (v \sc \partial w)_k = {i \over (2 \pi)^{d/2}}
\sum_{h \in \Zd} [v_{h} \sc \, (k-h)] w_{k - h}~; $$
this implies
\parn
\vbox{
\beq \la v \sc \partial w | w \ra_n =
\sum_{k \in \Zd} |k|^{2 n} \overline{(v \sc \partial w)_k} \, \sc w_k  \label{abimp} \feq
$$ =  - {i \over (2 \pi)^{d/2}}\sum_{h, k \in \Zd} |k|^{2 n} [\overline{v_{h}} \sc \, (k-h)]
(\overline{w_{k - h}} \sc w_k)~ =
- {i \over (2 \pi)^{d/2}}\sum_{h, \ell \in \Zd} |h + \ell|^{2 n} (\overline{v_{h}} \sc \, \ell)
(\overline{w_{\ell}} \sc w_{h + \ell})  $$
$$ = - {i \over (2 \pi)^{d/2}}\sum_{h \in V, \ell \in W, h + \ell \in W}
|h + \ell|^{2 n} (\overline{v_{h}} \sc \, \ell)
(\overline{w_{\ell}} \sc w_{h + \ell}) = - {i \over (2 \pi)^{d/2}} P_n((v)_k,(w)_k)~, $$
}
\noindent
which proves the thesis \rref{lappra}. In the above chain of equalities,
the third passage relies on a change of variable $k = h + \ell$, and the fourth passage
depends on the relations $v_h=0$ for $h \in \Zd \setminus V$,
$w_\ell=0$ for $\ell \in \Zd \setminus W$. \parn
\textsl{Step 3. Conclusion of the proof.}
We consider two nonzero families $(v_k)_{k \in V}$, $(w_k)_{k \in W}$ $\in \S$, and
define $v := \sum_{k \in V} v_k e_k$, $w := \sum_{k \in W} w_k e_k$. According to Steps
1 and 2, we have $\| v \|_n = N_n((v_k))$, $\| w \|_n = N_n((w_k))$,
$\la v \sc \partial w | w \ra_n =  (2 \pi)^{-d/2} P_n((v_k), (w_k))$; so, the inequality
$G_n \geqs  |\la v \sc \partial w | w \ra_n|/\| v \|_n \| w \|^2_n$ takes the form
(\ref{gnlowlow}-\ref{gnlow}). \fine
\vfill \eject \noindent
\appendix
\section{Some tools preparing the analysis of the function
$\boma{\GG_{n}}$}
\label{proveg}
In the sequel $d \in \{2,3,...\}$.
Let us fix some notations, to be used throughout the Appendices.
\begin{prop}
\textbf{Definition.}
\label{versk}
(i) $\teta : \reali \vain \{0,1\}$
is the Heaviside function such that $\teta(z) := 1$ if $z \in [0,+\infty)$ and $\teta(z) := 0$
if $z \in (-\infty,0)$. \parn
(ii) $\Gamma$ is the Euler Gamma function,
$\binom{\cdot}{\cdot}$ are the binomial coefficients. \parn
(iii) We put $\Sd := \{ u \in \reali^d~|~|u| = 1 \}$.
For each $p \in \reali^d \setminus \{0\}$, the versor of $p$ is
$\vers{p}  := \dd{p \over |p|} \in \Sd$.
\end{prop}
\begin{prop}
\label{lemmaf}
\textbf{Lemma.} For any function $f : \Zd_0 \vain \reali$ and $k \in \Zd_0$,
$\ro \in (1,+\infty)$, one has
\beq \sum_{h \in \Zd_{0 k}, |h| < \ro \op |k-h| < \ro} f(h)
= \sum_{h \in \Zd_{0 k},|h| < \ro} f(h) + \teta(|k-h| - \ro) f(k-h)~.
\label{tesif} \feq
\end{prop}
\textbf{Proof.} See \cite{cok}. \fine
\begin{prop}
\label{lemgs}
\textbf{Lemma.} For any $n \in (1,+\infty)$, the following holds. \parn
(i) Consider the function
\beq c_n : [0,4] \times [0,1] \vain [0,+\infty) \label{decn} \feq
$$ c_n(z,u) := \left\{ \barray{lll}
\dd{z (4 - z) [(1 - z u + z u^2)^{n/2} - (1 - u)^n  ]^2 \over
2 u^2 [ u^{2 n - 2} + (1 - u)^{2 n-2}]} & \mbox{if $u \in (0,1]$}~, \\
\dd{n^2 z (4 - z) (2 - z)^2 \over 8}~ & \mbox{if $u = 0$}~. \farray \right.
$$
This is well defined and continuous, which implies existence of
\beq C_n := \max_{z \in [0,4], u \in [0,1]} c_n(z,u) \in (0,+\infty)~. \label{decin} \feq
(ii) For all $p, q \in \reali^d$, one has
\beq |p \we q|^2 (|p+q|^n - |q|^{n})^2 \leqs {C_n \over 2} |p|^4 |q|^2
\Big[ |p|^{2 n - 2} + |q|^{2 n - 2} \Big]~. \label{onthesi} \feq
\end{prop}
\textbf{Proof.} (i) Well definedness and continuity of $c_n$ are
checked by elementary means, the main point being the computation of
$\lim_{u \vain 0^{+}} c_n(z,u)$.
\parn
(ii) Eq. \rref{onthesi} is obvious if $p=0$ or $q=0$, due to the vanishing of both sides;
hereafter we prove \rref{onthesi} for $p, q \in \reali^d \setminus \{0\}$.
Let $\tet(p,q) \equiv \tet \in [0,\pi]$ denote the convex angle between
$p$ and $q$; we have the relations
$$ |p \we q|^2 = |p|^2 |q|^2 \sin^2 \tet~, \qquad
|p + q|^2 = |p|^2 + |q|^2 + 2 |p| |q| \cos \tet~,$$
which imply
\beq {2 |p \we q|^2 (|p+q|^n - |q|^{n})^2 \over |p|^4 |q|^2
\Big[ |p|^{2 n - 2} + |q|^{2 n - 2} \Big]}
= {2 \sin^2 \tet [(|p|^2 + |q|^2 + 2 |p| |q| \cos \tet)^{n/2} - |q|^{n}]^2 \over
|p|^2 \Big[ |p|^{2 n - 2} + |q|^{2 n - 2} \Big]}~.  \label{eqgive} \feq
To go on, we define  $z \in [0,4]$, $u \in (0,1)$ through the equations
\beq \cos \tet = 1 - {z \over 2}~, \quad |p | = {u \over 1 - u} |q|~\feq
(note that $|p| = \ep |q|$ for a unique $\ep \in (0,+\infty)$;
on the other hand, the map $u \mapsto u/(1-u)$ is one-to-one between
$(0,1)$ and $(0,+\infty)$). Returning to \rref{eqgive}, after some computations we get
\beq {2 |p \we q|^2 (|p+q|^n - |q|^{n})^2 \over |p|^4 |q|^2
\Big[ |p|^{2 n - 2} + |q|^{2 n - 2} \Big]}
= c_n(z, u)~. \feq
But $c_n(z, u) \leqs C_n$, so we obtain the thesis
\rref{onthesi}. \fine
\begin{rema}
\textbf{Examples.} Let $c_n, C_n$ be defined as in the previous Lemma.
For $n=3,4,5,10$ we have the following numerical results, to be employed later:
\beq C_3= c_3(0.69603..., 0.46453...) = 14.814...\,;~\label{valc2} \feq
$$ C_4 = c_4(0.61987..., 0.47822...) = 58.460...\,; $$
$$ C_5 = c_5(0.55023..., 0.48569...) =  215.97...\,;~$$
$$ C_{10} = c_{10}(0.33289..., 0.49672...) = 1.3467... \times 10^5\,. $$
\end{rema}
\begin{prop}
\label{lemz}
\textbf{Lemma.} Let $\nu \in (d,+\infty)$. For any $\ro \in (2 \sqrt{d},+\infty)$, one has
\beq \sum_{h \in \Zd, |h| \geqs \ro} {1 \over |h|^{\nu}}~\leqs
{2 \pi^{d/2} \over \Gamma(d/2)}
\sum_{i=0}^{d-1} \left( \barray{c} d - 1 \\ i \farray \right) {d^{d/2-1/2-i/2} \over (\nu - i - 1)
(\ro - 2 \sqrt{d})^{\nu - i - 1}}~. \label{desnu}\feq
\end{prop}
\textbf{Proof.} This is just Lemma C.2 of  \cite{tre} (with the variable $\lan$ of the cited reference
related to $\ro$ by $\lan = \ro - 2 \sqrt{d}$). \fine
\begin{prop}
\textbf{Lemma.}
Let $\ro \in (1,+\infty)$ and $\varphi: [1,\ro) \vain \reali$.
Then, for each $k \in \reali^d$,
\beq
\sum_{h \in \Zd_0, |h| < \ro} (h \sc k)^2  \varphi(|h|) =
{|k|^2 \over d} \sum_{h \in \Zd_0, |h| < \ro} |h|^2 \varphi(|h|)~.
\label{idew} \feq
\end{prop}
\textbf{Proof.} See \cite{cok}. \fine
\begin{prop}
\label{deden}
\textbf{Definition.}
Let us introduce the domain
\beq \Do := \{ (c,\xi) \in \reali^2~|~c \in [-1,1], \, \xi \in [0,+\infty), \, (c,\xi) \neq (1,1) \}~;
\feq
furthermore, let $n \in \reali$. \parn
(i) We put
\beq D_n : \Do \vain [0,+\infty)~, \label{defn} \feq
$$ (c, \ep) \mapsto D_n(c,\ep) := \left\{ \barray{ll}
\dd{(1 - c^2)[1 - (1 - 2 c \ep + \ep^2)^{n/2}]^2 \over \ep^2 (1 - 2 c \ep + \ep^2)^n}~& \mbox{if
$\ep \neq 0$,} \\ n^2 (c^2 - c^4) & \mbox{if $\ep=0$}~; \farray \right. $$
\beq E_n : \Do \vain [0,+\infty)~,
\quad (c, \ep) \mapsto E_n(c,\ep) := {1 - c^2 \over (1 - 2 c \ep + \ep^2)^{n+1}}~. \label{den} \feq
($D_n$ is
$C^\infty$, as shown by an elementary analysis of the
term $\ep^{-2} [1 - (1 - 2 c \ep + \ep^2)^{n/2}]^2$;
$E_n$ already appeared in \cite{cok}, and is $C^\infty$ as well.) \parn
(ii) For $\ell = 0,1,2,...$, we put
\beq D_{n \ell}, E_{n \ell} :  [-1,1] \vain \reali,~D_{n \ell}(c) :=
{1 \over \ell!} {\partial^\ell D_n \over \partial \ep^\ell}(c,0),~
E_{n \ell}(c) :=
{1 \over \ell!} {\partial^\ell E_n \over \partial \ep^\ell}(c,0)~. \label{denel} \feq
(iii) For $t = 1,2,...$,
\beq Q_{n t}, R_{n t} :  \Do \vain \reali \feq
are the unique $C^\infty$ functions such that, for all $(c, \ep) \in \Do$,
\beq D_n(c, \ep) = \sum_{\ell=0}^{t-1} D_{n \ell}(c) \ep^\ell + Q_{n t}(c, \ep) \ep^t,
\quad
E_n(c, \ep) = \sum_{\ell=0}^{t-1} E_{n \ell}(c) \ep^\ell + R_{n t}(c, \ep) \ep^t.
\label{eqrt} \feq
(iv) For $t = 1,2,...$, we put
\beq \lam_{n t} := \min_{c \in [-1,1], \, \ep \in [0,1/2]} Q_{n t}(c, \ep)~,
\qquad \mu_{n t} := \min_{c \in [-1,1], \, \ep \in [0,1/2]} R_{n t}(c, \ep)~,
\label{munu} \feq
\beq \Lam_{n t} := \max_{c \in [-1,1], \, \ep \in [0,1/2]} Q_{n t}(c, \ep)~;
\qquad M_{n t} := \max_{c \in [-1,1], \, \ep \in [0,1/2]} R_{n t}(c, \ep)~. \label{mn} \feq
\end{prop}
\begin{rema}
\textbf{Remarks.} (i)
The first $D_{n \ell}$ functions are
\beq D_{n 0}(c) = n^2 (c^2 - c^4), \qquad D_{n 1}(c) = - n^2 c + (3 n^2 + n^3) c^3 - (2 n^2 + n^3) c^5,
 \label{fn02} \feq
$$ D_{n 2}(c)  := {n^2 \over 4} - ({13 \over 4} n^2 + {3 \over 2} n^3) c^2
+ ({20 \over 3} n^2 + {9 \over 2} n^3 + {7 \over 12} n^4) c^4
- ({11 \over 3} n^2 + 3 n^3 + {7 \over 12} n^4) c^6. $$
The first $E_{n \ell}$ functions are reported in \cite{cok}. \parn
(ii) In general, $D_{n \ell}$ and $E_{n \ell}$ are polynomials in $c$ of degrees $\ell + 4$
and $\ell + 2$, respectively; as functions of $c$, these have the same parity as $\ell$. \parn
(iii) Eq. \rref{eqrt} characterizes $Q_{n t}(c, \ep) \ep^t$ and $R_n(c, \ep) \ep^t$
as the reminders of two Taylor expansions. One can solve the equations in \rref{eqrt}
with respect to $Q_{n t}(c, \ep)$, $R_{n, t}(c, \ep)$; the expressions obtained in this way
can be used for the practical computation of these functions, and of their minima and maxima
defined by \rref{munu} \rref{mn}. Typically, the evaluation of the cited minima and maxima will be
numerical. \parn
(iv) For future use, we report here the minima and maxima, determined
numerically from the definitions \rref{munu} \rref{mn} with $n=3$, $t=8$ and
$n=4,5,10$, $t=6$:
\beq
\lam_{3 8} = -72.563...\,,~\Lam_{3 8} = 202.91...\,; \qquad
\mu_{3 8} = -159.61... \, ,~ M_{3 8} = 930.73... \, ; \label{m38} \feq
$$~\lam_{4 6} = -112.95... \, ,~
\Lam_{4 6} = 904.92... \, ; \quad \lam_{5 6} = -432.09... \, ,
~\Lam_{5 6} =  4970.4... \, ;~$$
$$ \lam_{10, 6} = -1.3678... \times 10^4... \, ,~ \Lam_{10, 6} = 5.0076... \times 10^6 \, . $$
(Some of the subsequent computations require as well the values of $m_{n 6}$, $M_{n 6}$
for $n=4,5,10$; these are reported in \cite{cok}.) \fine
\end{rema}
In the sequel we present a lemma
on a function of two vector variables $h,k$, to be used later
(see Eq.\rref{gmnd}); as indicated below,
this is related to the functions $D_n, E_n$ in \rref{den} and to their Taylor expansions.
\begin{prop}
\textbf{Lemma.} Let
$h, k \in \reali^d \setminus \{0\}$, $h \neq k$, and
let $\te(h,k) \equiv \te$ be the convex angle between them.
Furthermore, let $n \in \reali$;
then the following holds. \parn
(i) One has
\beq
|h \we k|^2 \left[{(|k|^n - |k-h|^n)^2 \over |h|^{2 n + 2} |k-h|^{2 n}} +
{(|k|^n - |h|^n)^2 \over |h|^{2 n} |k - h|^{2 n + 2} } \right] \label{fung}
\feq
$$ = {1 \over |h|^{2 n -2}}
\left[ D_n \Big(\cos \te, {|h| \over |k|} \Big)
+ \Big( 1 - {|h|^n \over |k|^n} \Big)^2 E_n \Big(\cos \te, {|h| \over |k|} \Big)
 \right]~. $$
(ii) Let $|k| \geqs 2 |h|$. For any $t \in \{1,2,...,\}$, Eq. \rref{fung} implies
\parn
\vbox{
$$ {1 \over |h|^{2 n -2}}
\left[
\sum_{\ell=0}^{t-1}
D_{n \ell} (\cos \te) {|h|^\ell \over |k|^\ell} +
\lam_{n t} {|h|^t \over |k|^t}
+ \Big( 1 - {|h|^n \over |k|^n} \Big)^2 \Big(\sum_{\ell=0}^{t-1}
E_{n \ell} (\cos \te) {|h|^\ell \over |k|^\ell} +
\mu_{n t} {|h|^t \over |k|^t} \Big) \right] $$
\beq
\leqs |h \we k|^2 \left[{(|k|^n - |k-h|^n)^2 \over |h|^{2 n + 2} |k-h|^{2 n}} +
{(|k|^n - |h|^n)^2 \over |h|^{2 n} |k - h|^{2 n + 2} } \right]
\label{eqtag} \feq
$$ \leqs {1 \over |h|^{2 n -2}}
\left[
\sum_{\ell=0}^{t-1}
D_{n \ell} (\cos \te) {|h|^\ell \over |k|^\ell} +
\Lam_{n t} {|h|^t \over |k|^t}
+ \Big( 1 - {|h|^n \over |k|^n} \Big)^2 \Big(\sum_{\ell=0}^{t-1}
E_{n \ell} (\cos \te) {|h|^\ell \over |k|^\ell} +
M_{n t} {|h|^t \over |k|^t} \Big)  \right] $$
}
\noindent
(note that $\cos \te = \vh \sc \vk$, with $\vers{~}$ denoting the versor).
\end{prop}
\textbf{Proof.} (i) We consider the function in the left hand side of \rref{fung}, and
reexpress it using the identities
\beq |h \we k|^2 = |h|^2 |k|^2 (1 - \cos^2 \te)~, \feq
$$ |k - h| = \sqrt{|k|^2 - 2 |k| |h| \cos \te + |h|^2} =
|k| \sqrt{1 - 2 \, \cos \te {|h| \over |k|}  + {|h|^2 \over |k|^2}}~, $$
$$ |k|^n - |h|^n = |k|^n \Big(1 - {|h|^n \over |k|^n}\Big)~; $$
these readily yield the thesis \rref{fung}. \parn
(ii) Eqs. \rref{eqrt} \rref{munu} \rref{mn} imply
\beq
\sum_{\ell=0}^{t-1} D_{n \ell}(c) \ep^\ell + \lam_{n t} \ep^t\leqs
D_n(c,\ep) \leqs \sum_{\ell=0}^{t-1} D_{n \ell}(c) \ep^\ell + \Lam_{n t} \ep^t~,
\label{enfor} \feq
$$ \sum_{\ell=0}^{t-1} E_{n \ell}(c) \ep^\ell + \mu_{n t} \ep^t\leqs
E_n(c,\ep) \leqs \sum_{\ell=0}^{t-1} E_{n \ell}(c) \ep^\ell + M_{n t} \ep^t
\quad \mbox{for $(c,\ep) \in [-1,1] \times [0,1/2]$}. $$
Let us apply these inequalities with $c := \cos \te$ and $\ep := |h|/|k|$
(noting that $0 \leqs \ep \leqs 1/2$, by the assumption $|k| \geqs 2 |h|$). In this
way, from Eqs. \rref{fung} and \rref{enfor} we readily get the thesis \rref{eqtag}.
\fine
To conclude, let us introduce some variants $\Dd_{n \ell}$ and
$\Ed_{n \ell}$ of the polynomials defined before ($\Ed_{n \ell}$
was already considered in \cite{cok}).
\begin{prop}
\label{dendel}
\textbf{Definition.} For $\ell = 0,2,... $, $\Dd_{n \ell d} \equiv \Dd_{n \ell}$
and $\Ed_{n \ell d} \equiv \Ed_{n \ell}$ are the polynomials obtained
from $D_{n \ell}$ and $E_{n \ell}$, replacing the term $c^2$ with $1/d$.
\end{prop}
\parn
\vbox{
\begin{rema}
\textbf{Example.} The expressions of $D_{n 0}$, $D_{n 2}$ in \rref{fn02} imply
\beq \Dd_{n 0}(c) = {n^2 \over d} - n^2 c^4~~, \label{en02d} \feq
$$ \Dd_{n 2}(c)  = {n^2 \over 4} - ({13 \over 4} n^2 + {3 \over 2} n^3) {1 \over d}
+ ({20 \over 3} n^2 + {9 \over 2} n^3 + {7 \over 12} n^4) c^4
- ({11 \over 3} n^2 + 3 n^3 + {7 \over 12} n^4) c^6~. $$
\end{rema}
}
\noindent
\section{The function $\boma{\GG_{n}}$}
\label{appeg}
Throughout the appendix $n \in ({d/2} + 1, + \infty)$.
For $k \in \Zd_0$, we recall the definition \rref{ggnd}
$$ \GG_{n}(k) :=
\sum_{h \in \Zd_{0 k}} {|h \we k|^2(|k|^n - |k-h|^n)^2
\over |h|^{2 n + 2} |k-h|^{2 n}} \in (0,+\infty)~, \qquad
(\Zd_{0 k} := \Zd \setminus \{0,k\})~.$$
\begin{prop}
\label{proggnd}
\textbf{Proposition.}
Let us choose a ''cutoff''
\beq \ro \in (2 \sqrt{d},+\infty)~; \label{cutoff} \feq
then, the following holds
(with the functions and quantities
$\Gg_{n}$, $\dG_{n}$,... mentioned in the sequel
depending parametrically on $d$ and $\ro$: $\Gg_{n}(k) \equiv \Gg_{n d}(k, \ro)$,
$\dG_{n} \equiv \dG_{n d}(\ro)$,...). \parn
(i) The function $\GG_{n}$ can be evaluated using the inequalities
\beq \Gg_{n}(k) < \GG_{n}(k) \leqs \Gg_{n}(k) + \dG_{n}
~~\mbox{for all $k \in \Zd_0$}~. \label{dgmnd} \feq
Here
\beq \Gg_{n}(k) :=
\sum_{h \in \Zd_{0 k}, |h| < \ro \op |k-h| < \ro} {|h \we k|^2(|k|^n - |k-h|^n)^2
\over |h|^{2 n + 2} |k-h|^{2 n}}~; \label{deagm} \feq
this function can be reexpressed as
\beq \Gg_{n}(k) = \hspace{-0.3cm}\sum_{h \in \Zd_0, |h | < \ro }\hspace{-0.3cm}
|h \we k|^2 \left[{(|k|^n - |k-h|^n)^2 \over |h|^{2 n + 2} |k-h|^{2 n}} + \teta(|k-h| - \ro)
{(|k|^n - |h|^n)^2 \over |h|^{2 n} |k - h|^{2 n + 2} } \right] \label{gmnd} \feq
(with $\teta$ as in Definition \ref{versk}).
If $|k| \geqs 2 \ro$, in Eq. \rref{gmnd} one can replace $\Zd_{0 k}$ with $\Zd_0$
and $\teta(|k-h| - \ro)$ with $1$.
Furthermore
\beq \dG_{n} :=
{2 \pi^{d/2} C_n \over \Gamma(d/2)}
\sum_{i=0}^{d-1} \Big( \barray{c} d - 1 \\ i \farray \Big) {d^{d/2-1/2-i/2} \over (2 n - 3 - i)
(\ro - 2 \sqrt{d})^{2 n - 3 - i}}~, \label{dedeg} \feq
with $C_n$ as in \rref{decin}. \parn
(ii) As in Remark \ref{rem34}, consider the reflection operators
$R_r$ ($r=1,...,d$)
and the permutation operators $P_{\si}$
($\si$ a permutation of $\{1,...,d\}$). Then
\beq \Gg_{n}(R_r k) = \Gg_{n}(k)~, \quad \Gg_{n}(P_\sigma k) = \Gg_{n}(k)~
\qquad \mbox{for each $k \in \Zd_0$} \label{claia} \feq
(so, the computation of $\Gg_{n}(k)$
can be reduced to the case $k_1 \geqs k_2 \geqs ... \geqs k_d \geqs 0$). \parn
(iii) Let $t \in \{2,4,...\}$. One has
\parn
\vbox{
$$
\sum_{\ell=0,2,...,t-2} {1 \over |k|^\ell} \left(\PP_{n \ell}(\vk) +
{\Pp_{n \ell}(\vk) \over |k|^{n}} +
{\Ps_{n \ell}(\vk) \over |k|^{2 n}}\right)
+ {1 \over |k|^t} \left(w_{n t} + {\wp_{n t} \over |k|^n} + {\ws_{n t} \over |k|^{2 n}} \right)
$$
\beq \leqs \Gg_{n}(k)  \label{takm} \feq
$$ {~} \hspace{-0.5cm} \leqs \hspace{-0.5cm}
\sum_{\ell=0,2,...,t-2} {1 \over |k|^\ell} \left(\PP_{n \ell}(\vk) +
{\Pp_{n \ell}(\vk) \over |k|^{n}} +
{\Ps_{n \ell}(\vk) \over |k|^{2 n}}\right)
+ {1 \over |k|^t} \left(W_{n t} + {\Wp_{n t} \over |k|^n} + {\Ws_{n t} \over |k|^{2 n}} \right)
~\mbox{for $k \in \Zd_0$, $|k| \geqs 2 \ro$}. $$
}
\noindent
In the above,
$\vk \in \Sd$
is the versor of $k$ (see Definition \ref{versk}). Furthermore,
\beq \PP_{n \ell}, \Pp_{n \ell}, \Ps_{n \ell} : \Sd \vain \reali,~\label{depnel} \feq
$$ \PP_{n \ell}(u) := \hspace{-0.2cm} \sum_{h \in \Zd_0, |h| < \ro} \hspace{-0.2cm}
{\Dd_{n \ell}(\vh \sc u) + \Ed_{n \ell}(\vh \sc u) \over |h|^{2 n - 2 - \ell}},~~
\Pp_{n \ell}(u) := - 2 \hspace{-0.4cm} \sum_{h \in \Zd_0, |h| < \ro} \hspace{-0.2cm}
{\Ed_{n \ell}(\vh \sc u)  \over |h|^{n - 2 - \ell}}~, $$
$$
\Ps_{n \ell}(u) := \hspace{-0.2cm} \sum_{h \in \Zd_0, |h| < \ro} \hspace{-0.2cm}
\Ed_{n \ell}(\vh \sc u) |h|^{2 + \ell}
\qquad (\mbox{$\Dd_{n \ell}$, $\Ed_{n \ell}$ as in Definition \ref{dendel}})~; $$
\beq
w_{n t} := (\lam_{n t} + \mu_{n t})
\hspace{-0.2cm} \sum_{h \in \Zd_0, |h| < \ro} \hspace{-0.2cm}
{1 \over |h|^{2 n - 2 - t}},~~
\wp_{n t} := - 2 \mu_{n t}
\hspace{-0.2cm} \sum_{h \in \Zd_0, |h| < \ro} \hspace{-0.2cm}
{1 \over |h|^{n - 2 - t}},~ \label{wnt} \feq
$$ \ws_{n t} := \mu_{n t}
\hspace{-0.2cm} \sum_{h \in \Zd_0, |h| < \ro} \hspace{-0.2cm} |h|^{2 + t}
\qquad (\mbox{$\lam_{n t}$, $\mu_{n t}$ as in Eq. \rref{munu}})~; $$
\beq
W_{n t} := (\Lam_{n t} + M_{n t})
\hspace{-0.2cm} \sum_{h \in \Zd_0, |h| < \ro} \hspace{-0.2cm}
{1 \over |h|^{2 n - 2 - t}},~~
\Wp_{n t} := - 2 M_{n t}
\hspace{-0.2cm} \sum_{h \in \Zd_0, |h| < \ro} \hspace{-0.2cm}
{1 \over |h|^{n - 2 - t}},~ \label{wwnt} \feq
$$ \Ws_{n t} := M_{n t}
\hspace{-0.4cm} \sum_{h \in \Zd_0, |h| < \ro} \hspace{-0.4cm} |h|^{2 + t}
\qquad (\mbox{$\Lam_{n t}$, $M_{n t}$ as in Eq. \rref{mn}})~. $$
For each $\ell$, $\PP_{n \ell}$, $\Pp_{n \ell}$ and $\Ps_{n \ell}$ are polynomial functions on $\Sd$; setting
\parn
\vbox{
\beq p_{n \ell} := \min_{u \in \Sd} \PP_{n \ell}(u), \quad
\pp_{n \ell} := \min_{u \in \Sd} \Pp_{n \ell}(u), \quad
\ps_{n \ell} := \min_{u \in \Sd} \Ps_{n \ell}(u), \label{mqnel} \feq
$$ P_{n \ell} := \max_{u \in \Sd} \PP_{n \ell}(u), \quad
\Pip_{n \ell} := \max_{u \in \Sd} \Pp_{n \ell}(u), \quad
\Pis_{n \ell} := \max_{u \in \Sd} \Ps_{n \ell}(u), $$
}
\noindent
one infers from \rref{takm} that
\parn
\vbox{
\beq
\sum_{\ell = 2,4,....,t-2} {1 \over |k|^\ell} \left(p_{n \ell} +
{\pp_{n \ell} \over |k|^{n}} +
{\ps_{n \ell} \over |k|^{2 n}}\right)
+ {1 \over |k|^t} \left(w_{n t} + {\wp_{n t} \over |k|^n} + {\ws_{n t} \over |k|^{2 n}} \right)
\leqs \Gg_{n}(k)  \label{tkm} \feq
$$ {~} \hspace{-0.5cm} \leqs \hspace{-0.5cm}
\sum_{\ell = 2,4,....,t-2} {1 \over |k|^\ell} \left(P_{n \ell} +
{\Pip_{n \ell} \over |k|^{n}} +
{\Pis_{n \ell}  \over |k|^{2 n}}\right)
+ {1 \over |k|^t} \left(W_{n t} + {\Wp_{n t} \over |k|^n} + {\Ws_{n t} \over |k|^{2 n}} \right)
~\mbox{for $k \in \Zd_0$, $|k| \geqs 2 \ro$}. $$
}
\noindent
Consider a sequence $(k_i)_{i=0,1,2,...}$ in $\Zd_0$; then the inequalities \rref{tkm}, with $t=2$, imply
\beq \Gg_n(k_i) \vain \PP_{n 0}(u) \qquad  \mbox{for $i \vain + \infty$, if $k_i \vain \infty$
and $\vers{k_i} \vain u \in \Sd$}~. \label{limu} \feq
Finally, we have
\beq \liminf_{k \in \Zd_0, k \vain \infty} \Gg_n(k) = p_{n 0}~,
\qquad \limsup_{k \in \Zd_0, k \vain \infty} \Gg_n(k) = P_{n 0}~. \label{limsup} \feq
(iv) Items (i) and (iii) imply
\beq \sup_{k \in \Zd_0} \Gg_{n}(k)  \leqs \sup_{k \in \Zd_0} \GG_{n}(k) \leqs
\Big(\sup_{k \in \Zd_0} \Gg_{n}(k) \Big) +
\dG_{n} < + \infty~. \label{impifi} \feq
\end{prop}
\textbf{Proof.} We fix a cutoff  $\ro$ as in \rref{cutoff}.
Our argument is divided in several steps; more precisely,
Steps 1-5 give proofs of statements (i)(ii), while
Steps 6-9 prove statements (iii)(iv).
The assumption \rref{cutoff} $\ro > 2 \sqrt{d}$ is essential in Step 3.
\parn
\textsl{Step 1. One has
\beq \GG_{n}(k)= \Gg_{n}(k) + \DG_{n}(k)~\quad \mbox{for all $k \in \Zd_0$}~, \label{decompg} \feq
where, as in \rref{deagm},
$\Gg_{n}(k) := \dd{\sum_{h \in \Zd_{0 k}, |h| < \ro \op |k-h| < \ro} {|h \we k|^2 (|k|^n - |k-h|^n)^2
\over |h|^{2 n + 2} |k-h|^{2 n}}}$,
while
\beq \DG_{n}(k) := \sum_{h \in \Zd_0, |h| \geqs \ro, |k-h| \geqs \ro}
 {|h \we k|^2 (|k|^n - |k-h|^n)^2
\over |h|^{2 n + 2} |k-h|^{2 n}}~
 \in (0,+\infty)~. \label{deagp} \feq}
The above decomposition follows noting that $\Zd_{0 k}$ is the disjoint union of
the domains of the sums defining $\Gg_n(k)$ and $\DG_{n}(k)$.
$\Gg_{n}(k)$ is finite, involving finitely many summands; $\DG_{n}(k)$ is
finite as well, since we know that $\GG_{n}(k) < +\infty$. \parn
\textsl{Step 2. For each $k \in \Zd_0$, one has the representation \rref{gmnd}
$$ \Gg_{n}(k)= \hspace{-0.3cm}\sum_{h \in \Zd_{0 k}, |h | < \ro }\hspace{-0.3cm}
|h \we k|^2 \left[{(|k|^n - |k-h|^n)^2 \over |h|^{2 n + 2} |k-h|^{2 n}} + \teta(|k-h| - \ro)
{(|k|^n - |h|^n)^2 \over |h|^{2 n} |k - h|^{2 n + 2} } \right]~. $$
If $|k| \geqs 2 \ro$, in the above
one can replace $\Zd_{0 k}$ with  $\Zd_{0}$
and $\teta(|k-h| - \ro)$ with $1$}. \parn
To prove \rref{gmnd} we reexpress the sum in Eq. \rref{deagm}, using
Eq. \rref{tesif} with $f(h) \equiv f_k(h) :=
\dd{|h \we k| (|k|^n - |k-h|^n)^2
\over |h|^{2 n + 2} |k-h|^{2 n}}$ (note that $f(k-h)$ contains
a term $|(k - h) \we k| = |h \we k|$). To go on, assume
$|k| \geqs 2 \ro$; then, for all $h \in \Zd_0$ with $|h| < \ro$ one has
$|k - h| \geqs |k| - |h| > \ro$;
this implies $h \neq k$ (i.e., $h \in \Zd_{0 k}$) and $\theta(|k-h| - \ro) = 1$, two facts which
justify the replacements indicated above.
\parn
\textsl{Step 3. For each $k \in \Zd_0$ one has
\beq 0 <  \DG_{n}(k) \leqs \dG_{n}~, \label{boundg} \feq
with $\dG_n$ as in Eq. \rref{dedeg}.}
The obvious relation $0 < \DG_{n}(k)$ was already noted; in the sequel we prove that
$\DG_{n}(k) \leqs \dG_{n}$.
The definition \rref{deagp} of $\DG_{n}(k)$ contains the term
$|h \we k|^2 (|k|^n - |k-h|^n)^2$, for which we have:
\parn
\vbox{
\beq |h \we k|^2 (|k|^n - |k-h|^n)^2 = |h \we (k-h)|^2 (|k|^n - |k-h|^n)^2  \label{ontheses} \feq
$$ \leqs {C_n \over 2} |h|^4 |k-h|^2
\Big[ |h|^{2 n - 2} + |k-h|^{2 n - 2} \Big] $$
}
\noindent
(the last inequality follows from \rref{onthesi}, with $p=h$ and $q=k-h$).
Inserting \rref{ontheses} into \rref{deagp}, we obtain
\beq \DG_{n}(k) \leqs {C_n \over 2} \sum_{h \in \Zd_0, |h| \geqs \ro, |k-h| \geqs \ro}
{|h|^{2 n - 2} + |k-h|^{2 n - 2}
\over |h|^{2 n - 2} |k-h|^{2 n - 2}}  \feq
$$ = {C_n \over 2} \, \Big( \sum_{h \in \Zd_0, |h| \geqs \ro, |k-h| \geqs \ro} {1 \over |k-h|^{2 n - 2}} +
\sum_{h \in \Zd_0, |h| \geqs \ro, |k-h| \geqs \ro} {1 \over |h|^{2 n - 2}} \Big)~. $$
The domain of the above two sums is contained in each one of the sets
$\{ h \in \Zd~|~|h | \geqs \ro \}$ and $\{ h \in \Zd~|~|k - h | \geqs \ro \}$; so,
\beq \DG_{n}(k,\ro) \leqs
{C_n \over 2} \, \Big(\sum_{h \in \Zd_0, |k-h| \geqs \ro} {1 \over |k-h|^{2 n - 2}} +
\sum_{h \in \Zd_0, |h | \geqs \ro} {1 \over |h|^{2 n - 2}} \Big)~. \feq
Now, the change of variable $h \mapsto k - h$ in the first sum shows that it is equal to
the second one, so
\beq \DG_{n}(k) \leqs C_n \sum_{h \in \Zd_0, |h | \geqs \ro}
{1 \over |h|^{2 n - 2}}~. \label{eefg} \feq
Finally, Eq. \rref{eefg} and Eq. \rref{desnu} with $\nu = 2 n - 2$ give
$$ \DG_{n}(k)
\leqs
{2 \pi^{d/2} C_n \over \Gamma(d/2)}
\sum_{i=0}^{d-1} \left( \barray{c} d - 1 \\ i \farray \right) {d^{d/2-1/2-i/2} \over (2 n - 3 -i)
(\ro - 2 \sqrt{d})^{2 n - 3 -i}}~=\dG_{n}~\mbox{as in \rref{dedeg}}~. $$
\parn
\textsl{Step 4. One has the inequalities \rref{dgmnd}
$\Gg_{n}(k) < \GG_{n}(k) \leqs \Gg_{n}(k) + \dG_{n}$}.
These relations follow immediately from the decomposition
\rref{decompg}
$\GG_{n}(k)= \Gg_{n}(k) + \DG_{n}(k)$
and from the bounds \rref{boundg} on $\DG_{n}(k)$. \parn
\parn
\textsl{Step 5. One has the equalities \rref{claia}
$\Gg_{n}(R_r k) = \Gg_{n}(k)$, $\Gg_{n}(P_\sigma k) = \Gg_{n}(k)$,
involving the reflection and permutation operators  $R_r, P_\si$.} Again,
we can invoke the argument employed for the analogous properties
of the function $\KK_{n}$ in \cite{cok}. \parn
\textsl{Step 6. Let $t \in \{2,4,...\}$. One has the inequalities \rref{takm} for
$\Gg_n$.}
As an example, for any $k \in \Zd_0$ with $|k| \geqs 2 \ro$ we prove the upper bound \rref{takm}
$$\Gg_{n}(k) \leqs \hspace{-0.5cm}
\sum_{\ell=0,2,...,t-2} {1 \over |k|^\ell} \left(\PP_{n \ell}(\vk) +
{\Pp_{n \ell}(\vk) \over |k|^{n}} +
{\Ps_{n \ell}(\vk) \over |k|^{2 n}}\right)
+ {1 \over |k|^t} \left(W_{n t} + {\Wp_{n t} \over |k|^n} + {\Ws_{n t} \over |k|^{2 n}} \right). $$
Since $|k| \geqs 2 \ro$, we can express $\Gg_n(k)$
via Eq. \rref{gmnd}, replacing therein $\Zd_{0 k}$ with $\Zd_0$ and
$\theta(|k-h| - \ro)$ with $1$ (see the final statement in Step 2). So,
\beq \Gg_{n}(k)= \hspace{-0.3cm}\sum_{h \in \Zd_{0}, |h | < \ro }\hspace{-0.3cm}
|h \we k|^2 \left[{(|k|^n - |k-h|^n)^2 \over |h|^{2 n + 2} |k-h|^{2 n}}
+  {(|k|^n - |h|^n)^2 \over |h|^{2 n} |k - h|^{2 n + 2} } \right]~. \feq
In this expression we insert the upper bound of Eq. \rref{eqtag},
writing therein $\cos \te = \vh \sc \vk$
(note that \rref{eqtag} can be used, since $|h|/|k| < \ro/(2 \ro) < 1/2$ for each $h$ in the sum).
After some elementary manipulations (such as expanding
the square $(1 - {|h|^n / |k|^n} )^2$, and reorganizing the terms
that arise in this way), we conclude
$$ \Gg_{n}(k) \leqs
\sum_{\ell=0,1,...,t-1} {1 \over |k|^\ell} \left(\PP_{n \ell}(\vk) +
{\Pp_{n \ell}(\vk) \over |k|^{n}} +
{\Ps_{n \ell}(\vk) \over |k|^{2 n}}\right)
+ {1 \over |k|^t} \left(W_{n t} + {\Wp_{n t} \over |k|^n} + {\Ws_{n t} \over |k|^{2 n}} \right), $$
where $W_{n t}, \Wp_{n t}, \Ws_{n t}$ are as in \rref{wwnt} and, for each $\ell \in \{0,...,t-1\}$,
we have provisionally defined
\beq \PP_{n \ell}, \Pp_{n \ell}, \Ps_{n \ell} : \Sd \vain \reali,~\label{deppnel} \feq
$$ \PP_{n \ell}(u) := \hspace{-0.4cm} \sum_{h \in \Zd_0, |h| < \ro} \hspace{-0.4cm}
{D_{n \ell}(\vh \sc u) + E_{n \ell}(\vh \sc u) \over |h|^{2 n - 2 - \ell}},~~
\Pp_{n \ell}(u) := - 2 \hspace{-0.4cm} \sum_{h \in \Zd_0, |h| < \ro} \hspace{-0.4cm}
{E_{n \ell}(\vh \sc u)  \over |h|^{n - 2 - \ell}}~, $$
$$
\Ps_{n \ell}(u) := \hspace{-0.4cm} \sum_{h \in \Zd_0, |h| < \ro} \hspace{-0.4cm}
E_{n \ell}(\vh \sc u) |h|^{2 + \ell}
\qquad (\mbox{$E_{n \ell}$, $D_{n \ell}$ as in Definition \ref{deden}})~. $$
Now, the thesis follows if we prove the following relations:
\beq \PP_{n \ell}(u) = 0, ~~\Pp_{n \ell}(u)=0,~~ \Ps_{n \ell}(u)=0 \quad
~\mbox{for $\ell \in \{1,3,...,t-1\}$, $u \in \Sd$}~; \label{b35} \feq
\beq {~} \hspace{-1cm}~\PP_{n \ell}(u), \Pp_{n \ell}(u), \Ps_{n \ell}(u)~\mbox{are as in \rref{depnel},
for $\ell \in \{0,2,4,...,t-2\}$, $u \in \Sd$}~.
\label{b36}\feq
The relations \rref{b35} are proved recalling
that, for $\ell$ odd, the functions $c \mapsto E_{n \ell}(c)$, $D_{n \ell}(c)$ are odd as well; this implies
that the general term of the sum \rref{deppnel} changes its sign under
a transformation $h \mapsto -h$. \parn
Now, let us prove \rref{b36} for any even $\ell$. As an example, we consider
the case of $\PP_{n \ell}$; the sum defining it in \rref{deppnel} contains
the even polynomials
\beq
D_{n \ell}(c) = \sum_{j=0,2,...,\ell+4} D_{n \ell j} c^j,
\qquad
E_{n \ell}(c) = \sum_{j=0,2,...,\ell+2} E_{n \ell j} c^j, \feq
so \rref{deppnel} implies
\beq \PP_{n \ell}(u) =
 \sum_{j=0,2,...,\ell+4} \!\!\! D_{n \ell j} \hspace{-0.3cm} \sum_{h \in \Zd_0, |h | < \ro }
{(\vh \sc u)^j \over |h|^{2 n - 2 - \ell} } +
\!\!\!  \sum_{j=0,2,...,\ell+2} E_{n \ell j} \hspace{-0.3cm} \sum_{h \in \Zd_0, |h | < \ro }
{(\vh \sc u)^j \over |h|^{2 n - 2 - \ell} } \, ;  \label{b40} \feq
in particular, for the $j=2$ terms in both sums above we have (writing $\vh = h/|h|$)
\beq \sum_{h \in \Zd_0, |h | < \ro } {(\vh \sc u)^2 \over |h|^{2 n - 2 - \ell} }~=
\sum_{h \in \Zd_0, |h | < \ro } {(h \sc u)^2 \over |h|^{2 n - \ell} }~=
{1 \over d} \sum_{h \in \Zd_0, |h| < \ro }
{1 \over |h|^{2 n - 2 - \ell} }~, \label{b41} \feq
where the last passage follows from the identity \rref{idew} (with $k$ replaced by $u$
and $\varphi(|h|) = 1/|h|^{2 n - \ell}$). Eqs. \rref{b40} \rref{b41} imply
\beq \PP_{n \ell}(u) =
\sum_{j=0,4,6,...,\ell+4} D_{n \ell j} \sum_{h \in \Zd_0, |h | < \ro }
{(\vh \sc u)^j \over |h|^{2 n - 2 - \ell} }
+ {D_{n \ell 2} \over d} \sum_{h \in \Zd_0, |h | < \ro } {1 \over |h|^{2 n - 2 - \ell} }
\label{dacomp} \feq
$$ + \sum_{j=0,4,6...,\ell+2} E_{n \ell j} \sum_{h \in \Zd_0, |h | < \ro }
{(\vh \sc u)^j \over |h|^{2 n - 2 - \ell} }
+ {E_{n \ell 2} \over d} \sum_{h \in \Zd_0, |h | < \ro } {1 \over |h|^{2 n - 2 - \ell} }~. $$
On the other hand, Definition \ref{dendel} of $\Dd_{n \ell}$, $\Ed_{n \ell}$ prescribes
\beq
\Dd_{n \ell}(c) = \!\!\! \sum_{j=0,4,6,...,\ell+4} D_{n \ell j} c^j + {D_{n \ell 2} \over d}~,
\qquad
\Ed_{n \ell}(c) = \!\!\! \sum_{j=0,4,6,...,\ell+2} E_{n \ell j} c^j + {E_{n \ell 2} \over d}~; \feq
comparing this with \rref{dacomp}, we conclude
$$ \PP_{n \ell}(u) = \hspace{-0.4cm} \sum_{h \in \Zd_0, |h| < \ro} \hspace{-0.4cm}
{\Dd_{n \ell}(\vh \sc u|) + \Ed_{n \ell}(\vh \sc u) \over |h|^{2 n - 2 - \ell}},~~ \mbox{as in \rref{depnel}}.
$$
So, statement \rref{b36} is proved for $\PP_{n \ell}$; one proceeds similarly for $\Pp_{n \ell}$
and $\Ps_{n \ell}$. \parn
\textsl{Step 7. Let $t \in \{2,4,...\}$. For $\ell \in \{0,2,4,...,t-2\}$,
$\PP_{n \ell}$, $\Pp_{n \ell}$ and $\Ps_{n \ell}$ are polynomial function on $\Sd$;
considering their minima and maxima
$p_{n \ell}$, $P_{n \ell}$, etc., one infers from \rref{takm} the inequalities
\rref{tkm}
$$
\sum_{\ell=0,2,4,....,t-2} {1 \over |k|^\ell} \left(p_{n \ell} +
{\pp_{n \ell} \over |k|^{n}} +
{\ps_{n \ell} \over |k|^{2 n}}\right)
+ {1 \over |k|^t} \left(w_{n t} + {\wp_{n t} \over |k|^n} + {\ws_{n t} \over |k|^{2 n}} \right)
\leqs \Gg_{n}(k) $$
$$ {~} \hspace{-0.5cm} \leqs \hspace{-0.5cm}
\sum_{\ell=0,2,4,....,t-2} {1 \over |k|^\ell} \left(P_{n \ell} +
{\Pip_{n \ell} \over |k|^{n}} +
{\Pis_{n \ell}  \over |k|^{2 n}}\right)
+ {1 \over |k|^t} \left(W_{n t} + {\Wp_{n t} \over |k|^n} + {\Ws_{n t} \over |k|^{2 n}} \right)
\quad \mbox{for $|k| \geqs 2 \rho$}. $$}
The polynomial nature of the functions $\PP_{n \ell}$, $\Pp_{n \ell}$ and
$\Ps_{n \ell}$ follows
from their definition \rref{depnel} in terms of the polynomials $\Ed_{n \ell}, \Dd_{n \ell}$.
The inequalities \rref{tkm} are obvious. \parn
\textsl{Step 8.
Consider a sequence $(k_i)_{i=0,1,2,...}$ in $\Zd_0$; then the inequalities \rref{tkm}, with $t=2$, imply
statement \rref{limu}
$$ \Gg_n(k_i) \vain \PP_{n 0}(u) \qquad  \mbox{for $i \vain + \infty$, if $k_i \vain \infty$
and $\vers{k_i} \vain u \in \Sd$}~. $$
Finally, we have the results \rref{limsup}
$$ \liminf_{k \in \Zd_0, k \vain \infty} \Gg_n(k) = p_{n 0}~,
\qquad \limsup_{k \in \Zd_0, k \vain \infty} \Gg_n(k) = P_{n 0}~. $$}
To prove all this we start from any sequence $(k_i)_{i=0,1,...}$ in $\Zd_0$ and
note that \rref{takm},
with $t=2$ and $k = k_i$, gives
\beq \PP_{n 0}(\vers{k_i}) +
{\Pp_{n 0}(\vers{k_i}) \over |k_i|^{n}} +
{\Ps_{n 0}(\vers{k_i}) \over |k_i|^{2 n}}
+ {1 \over |k_i|^2} \left(w_{n 2} + {\wp_{n 2} \over |k_i|^n} + {\ws_{n 2} \over |k_i|^{2 n}} \right)
\leqs \Gg_{n}(k_i) \label{takm2} \feq
$$ \leqs \PP_{n 0}(\vers{k_i}) +
{\Pp_{n 0}(\vers{k_i}) \over |k_i|^{n}} +
{\Ps_{n 0}(\vers{k_i}) \over |k_i|^{2 n}}
+ {1 \over |k_i|^2} \left(W_{n 2} + {\Wp_{n 2} \over |k_i|^n} + {\Ws_{n 2} \over |k_i|^{2 n}} \right)
\mbox{for $|k_i| \geqs 2 \ro$}~. $$
Now, assume $k_i \vain \infty$ and $\vers{k_i} \vain u \in \Sd$; then,
both the lower and the upper bounds to $\Gg_n(k_i)$ in \rref{takm2} tend to
$\PP_{n 0}(u)$ and we obtain Eq. \rref{limu}. \parn
Let us pass to the proof of Eq. \rref{limsup};
as an example, we derive the statement about $\limsup_{k \vain \infty} \Gg_n(k)$.
By definition,
\beq \limsup_{k \in \Zd_0, k \vain \infty} \Gg_n(k) = \sup_{(k_i) \in \Cgot} \, \lim_{i \vain +\infty}
\Gg_n(k_i)~, \feq
$$ \Cgot := \{ \mbox{sequences $(k_i)_{i =0,1,2,...}$ in $\Zd_0$ such that $k_i \vain \infty$,
$\lim_{i \vain +\infty} \Gg_n(k_i)$ exists} \, \}~. $$
Consider any sequence $(k_i) \in \Cgot$; applying the upper bound in
Eq. \rref{tkm}, with $t=2$ and $k=k_i$, we get
\beq \Gg_{n}(k_i) \leqs P_{n 0} +
{\Pip_{n 0} \over |k_i|^{n}} +
{\Pis_{n 0} \over |k_i|^{2 n}}
+ {1 \over |k_i|^2} \left(W_{n 2} + {\Wp_{n 2} \over |k_i|^n} + {\Ws_{n 2} \over |k_i|^{2 n}} \right) \feq
for all $i$ such that $|k_i| \geqs 2 \ro$. Let $i \vain +\infty$; then $k_i \vain \infty$, and
the previous inequality implies
\beq \lim_{i \vain + \infty} \Gg_{n}(k_i) \leqs P_{n 0}~. \label{limsupa} \feq
Now, let $\up \in \Sd$ be such that
\beq \PP_{n \ell}(\up) = P_{n 0}~, \label{upsu} \feq
and let us consider a sequence $(k_i)_{i=0,1,2,...}$ in $\Zd_0$ such that
\beq k_i \vain \infty~,~\vers{k_i} \vain \up \qquad \mbox{for $i \vain + \infty$} \label{sucht} \feq
(e.g., $k_i := ([i \up_{1}],...,[i \up_{d}])$, where
$[~]$ is the integer part). Eqs. \rref{sucht} \rref{limu} and \rref{upsu} give
\beq \lim_{i \vain + \infty} \Gg_{n}(k_i) = P_{n 0}~. \label{limsupb} \feq
The results \rref{limsupa} and \rref{limsupb} imply
$\limsup_{k \vain \infty} \Gg_n(k) \leqs P_{n 0}$ and $\limsup_{k \vain \infty} \Gg_n(k) \geqs P_{n 0}$,
respectively, yielding the desired relation
\beq \limsup_{k \in \Zd_0, k \vain \infty} \Gg_n(k) = P_{n 0}~. \feq
\textsl{Step 9. Proof of the inequalities \rref{impifi}}
$$  \sup_{k \in \Zd_0} \Gg_{n}(k)
\leqs \sup_{k \in \Zd_0} \GG_{n}(k) \leqs \Big(\sup_{k \in \Zd_0} \Gg_{n}(k) \Big) +
\dG_{n} < + \infty~. $$
The first two inequalities are obvious consequences of the relations
\rref{dgmnd} $\Gg_{n}(k) < \GG_{n}(k) \leqs \Gg_{n}(k) + \dG_{n}$; the third inequality above
holds if we show that
\beq \sup_{k \in \Zd_0} \Gg_{n}(k) < + \infty~, \feq
and this follows from the finiteness of $\limsup_{k \vain \infty} \Gg_n(k)$
(see Step 8).
\fine
\vfill \eject \noindent
\section{Appendix. The upper bounds $\boma{G^{+}_{n}}$, for
$\boma{d=3}$ and $\boma{n=3,4,5,10}$}
\label{appe345g}
\vskip 2cm \noindent
Eq. \rref{desup} defines $G^{+}_n$ in terms
of $\sup_{k \in \Zt_0} \GG_n(k)$, or of any upper approximant for this sup.
In all the cases analyzed hereafter, we produce both an upper and a lower
approximant; the lower one is given only to indicate the
uncertainty in our evaluation of $\sup \GG_n$.
\vskip 0.5cm \noindent
\textbf{Some details on the
evaluation of $\boma{\GG_{3}}$ and of its sup.}
Among the examples
presented here, the case of $\GG_{3}$ is
the one requiring more expensive computations. \parn
To evaluate $\GG_{3}$, we apply Proposition \ref{proggnd} with a fairly large cutoff
\beq \ro = 20~;  \feq
thus, following Proposition \ref{proggnd}, we must often sum over the set
$\{h \in \Zt_0~|~|h| < 20 \}$. \parn
Eq. \rref{dedeg} (with the value of $C_3$
in \rref{valc2}) gives
\beq \dG_{3} = 12.478...~, \label{dg3} \feq
and it remains to evaluate the function $\Gg_{3}$. \parn
To compute  $\Gg_{3}(k)$, we start from the $k's$ in $\Zt_0$ with $|k| < 2 \ro = 40$.
Using directly the definition \rref{gmnd}
for all such $k$'s ({\footnote{In fact, due to the symmetry properties \rref{claia},
computation of $\Gg_3(k)$ can be limited to points $k$
such that $k_1 \geqs k_2 \geqs k_3 \geqs 0$.}}), we obtain
\beq \max_{k \in \Zt_0, |k| < 40} \Gg_{3}(k) = \Gg_{3}(9,9,9) = 34.901...~.  \label{c3} \feq
Let us pass to the case $|k| \geqs 40$. Here, our main tool is the upper bound in \rref{tkm} with $t=8$;
after some computations, this gives
\vfill \eject \noindent
{~}
\vskip -2cm
\parn
\vbox{
\beq \Gg_3(k) \leqs
33.725 + {1070.6 \over |k|^2} - {3337.9 \over |k|^3}
+ {2.9764 \times 10^5 \over |k|^4} - {2.6596 \times 10^6 \over |k|^5} \label{e3} \feq
$$  +
{1.3451 \times 10^8 \over |k|^6} - {1.7663 \times 10^9 \over |k|^7}
+ {2.5858 \times 10^{12} \over |k|^8} - {1.0476 \times 10^{12} \over |k|^9}
+ {4.7461 \times 10^{12} \over |k|^{10}}  $$
$$ - {2.3621 \times 10^{16} \over |k|^{11}}
+ {3.1212 \times 10^{15} \over |k|^{12}} + {7.2378 \times 10^{19} \over |k|^{14}}
\leqs 34.792 \quad \mbox{for $k \in \Zt_0$, $|k| \geqs 40$} $$
}
({\footnote{Let us give some supplementary information on the computations
yielding \rref{e3}. The $t=8$ upper bound in Eq. \rref{tkm} reads:
$$ \Gg_3(k) \leqs \hspace{-0.2cm}
\sum_{\ell=0,2,4,6} {1 \over |k|^\ell} \left(P_{3 \ell} +
{\Pip_{3 \ell} \over |k|^{3}} +
{\Pis_{3 \ell}  \over |k|^{6}}\right)
+ {1 \over |k|^8} \left(W_{3 8} + {\Wp_{3 8} \over |k|^3} + {\Ws_{3 8} \over |k|^{6}} \right)
~\mbox{for $k \in \Zd_0$, $|k| \geqs 40$}. $$
The constants $W_{3 8}$, $\Wp_{3 8}$, $\Ws_{3 8}$ are computed directly
from the definition \rref{wwnt}
(this requires previous knowledge of $M_{3 8} = 930.73...$ and $\Lam_{3 8} = 202.91...$,
see Eq. \rref{m38}). For $\ell=0,2,4,6$,
$P_{3 \ell}$, $\Pip_{3 \ell}$ and $\Pis_{3 \ell}$
are the maxima of the polynomial functions
$\PP_{3 \ell}$, $\Pp_{3 \ell}$ and $\Ps_{3 \ell}$ on $\St$;
for example, Eq. \rref{depnel} with $n=3$, $\ell=0$
gives
$$ \PP_{3 0}(u) = 58.311...
- 39.076... \,(u_1^2 u_2^2 + u_1^2 u_3^2 + u_2^2 u_3^2)
- 34.683... \,(u_1^4 + u_2^4 + u_3^4) $$
for all $u \in \St$, and one finds that $P_{3 0} =
\PP_{3 0}(1/\sqrt{3}, 1/\sqrt{3}, 1/\sqrt{3}) = 33.724...$.
Computing the other polynomials mentioned above and their maxima,
and rounding up from above the numerical outputs, we obtain
the first inequality \rref{e3} $\Gg_3(k) \leqs
33.725 + 1070.6 \, |k|^{-2} + $ etc.,
holding for $|k| \geqs 40$; on the other hand,
$33.725 + 1070.6\, |k|^{-2} + $ etc. $\leqs 34.792$
for all such $k$'s, which explains the second inequality \rref{e3}.}}).
(For completeness, we mention that the $t=8$ lower bound in
\rref{tkm} and Eq. \rref{limsup} imply
$\inf_{k \in \Zt_0, |k| \geqs 40} \Gg_{3}(k) = \liminf_{k \in \Zt_0,
k \vain \infty} \Gg_{3}(k) = 23.627...$, while $\limsup_{k \in \Zt_0,
k \vain \infty} \Gg_{3}(k) = 33.724...\,.$
({\footnote{Let us explain how to derive these statements. First of all,
Eq. \rref{limsup} gives
$$ \liminf_{k \in \Zd_0, k \vain \infty} \Gg_3(k) = p_{3 0}~,
\qquad \limsup_{k \in \Zd_0, k \vain \infty} \Gg_3(k) = P_{3 0}~,  $$
where $p_{3 0}$ and $P_{3 0}$ are the minimum and the maximum of the polynomial
$\PP_{3 0}$ over $\St$. The explicit expression of $\PP_{3 0}$ is given
in the previous footnote; it turns out that
$p_{3 0} = \PP_{3 0}(1,0,0) = 23.627...$ and (as stated before)
$P_{3 0} = \PP_{3 0}(1/\sqrt{3}, 1/\sqrt{3}, 1/\sqrt{3}) = 33.724...~$. \parn
Now, let us use the lower bound \rref{tkm} with $n=3$, $t=8$; computing
all the necessary constants, after some round up we get
$$ \Gg_3(k) \geqs p_{3 0} + {1042.9 \over |k|^2} - {3338.0 \over |k|^3}
+ {2.9617 \times 10^5 \over |k|^4}  - {2.6755 \times 10^6 \over |k|^5}
+ {1.3449 \times 10^8 \over |k|^6} - {1.7822 \times 10^9 \over |k|^7} $$
$$
- {5.2231 \times 10^{11} \over |k|^8} - {1.0729 \times 10^{12} \over |k|^9}
+ {4.6822 \times 10^{12} \over |k|^{10}} + {4.0510 \times 10^{15} \over |k|^{11}}
+ {3.0213 \times 10^{15} \over |k|^{12}} - {1.2413 \times 10^{19} \over |k|^{14}} $$
for $k \in \Zt_0$, $|k| \geqs 40$.
On the other hand, one has
$1042.9 \, |k|^{-2} - 3338.0 |k|^{-3} + ... \geqs 0$
for $|k| \geqs 40$; so, $\inf_{k \in \Zt_0, |k| \geqs 40} \Gg_3(k) \geqs p_{3 0}$.
It is obvious that $\inf_{k \in \Zt_0, |k| \geqs 40} \Gg_3(k) \leqs
\liminf_{k \in \Zt_0, k \vain \infty} \Gg_3(k)$; the latter equals $p_{3 0}$,
thus $\inf \Gg_3 = \liminf \Gg_3 = p_{3 0}$.}}).) \parn
The results \rref{c3} \rref{e3} yield
\beq \sup_{k \in \Zt_0} \Gg_{3}(k) = \Gg_{3}(9,9,9) = 34.901...~. \label{summag} \feq
We now pass to the function $\GG_{3}$; according to \rref{impifi} we have
$\sup_{k \in \Zt_0} \Gg_{3 }(k) \leqs
\sup_{k \in \Zt_0} \GG_{3}(k) \leqs \Big(\sup_{k \in \Zt_0} \Gg_{3}(k) \Big) +
\dG_{3}$, and the numerical results \rref{dg3} \rref{summag} give
\beq 34.901 < \sup_{k \in \Zt_0} \GG_{3}(k) < 47.381~. \label{sugit} \feq
(The uncertainty on this sup is fairly large, due to the
value of $\dG_3$ in \rref{dg3}; the error $\dG_3$
could be significantly reduced choosing a cutoff $\ro \gg 20$,
but the related computations would be much more expensive.) \parn
\textbf{The upper bound $\boma{G^{+}_{3}}$.}
According to the definition \rref{desup}, we have
\beq G^{+}_{3} = {1 \over (2 \pi)^{3/2}}
\sqrt{\sup_{k \in \Zt_0} \GG_{3}(k)} \quad \mbox{(or any upper approximant for this)}~. \feq
Due to \rref{sugit}, we can take $G^{+}_3 = (2 \pi)^{-3/2} \sqrt{47.381}\,$; rounding up to  three digits we can write
\beq G^{+}_{3} =  0.438 \, , \feq
as reported in \rref{boug}. \salto
\textbf{Preparing the examples with $\boma{n=4,5,10}$.}
To evaluate $\GG_n$ for the cited values of $n$, we apply Proposition \ref{proggnd} with a cutoff
\beq \ro = 10~;  \feq
thus, all sums over $h$ in Proposition \ref{proggnd} are over the set
$\{ h \in \Zt_0~|~|h| < 10 \}$.
\salto
\textbf{Some details on the
evaluation of $\boma{\GG_{4}}$ and of its sup.}
Eq. \rref{dedeg} (with the value of $C_4$
in \rref{valc2}) gives
\beq \dG_{4} = 1.2626...~, \label{dg4} \feq
and it remains to evaluate the function $\Gg_{4}$.
\parn
To compute  $\Gg_{4}(k)$, we start from the $k's$ in $\Zt_0$ with $|k| < 2 \ro = 20$.
Using directly the definition \rref{gmnd}
for all such $k$'s, we obtain
\beq \max_{k \in \Zt_0, |k| < 20} \Gg_{4}(k) = \Gg_{4}(2,1,0) = 56.628...~.  \label{c4} \feq
Let us pass to the case $|k| \geqs 20$. Here we use the upper bound in \rref{tkm} with $t=6$,
giving
\parn
\vbox{
\beq \Gg_4(k) \leqs 31.379 + {193.19 \over |k|^2}
+ {3740.3 \times 10^5 \over |k|^4} + {1.1291 \times 10^7 \over |k|^6}
- {8.6865 \times 10^{6} \over |k|^8} \label{e4} \feq
$$
- {6.3946 \times 10^{10} \over |k|^{10}}
+ {2.4366 \times 10^{5} \over |k|^{12}} + {2.0079 \times 10^{14} \over |k|^{14}}
\leqs 32.056 \quad \mbox{for $k \in \Zt_0$, $|k| \geqs 20$}~. $$
}
\noindent
(For completeness we mention that the $t=6$ lower bound in
\rref{tkm} and Eq. \rref{limsup} imply
$\inf_{k \in \Zt_0, |k| \geqs 20} \Gg_{4}(k) =
\liminf_{k \in \Zt_0, k \vain \infty} \GG_4(k) = 11.716...$,
$\limsup_{k \in \Zt_0, k \vain \infty} \Gg_{4}(k) = 31.378...$).
\parn
The results \rref{c4} \rref{e4} yield
\beq \sup_{k \in \Zt_0} \Gg_{4}(k) = \Gg_{4}(2,1,0) = 56.628...~.   \label{summag4} \feq
We now pass to the function $\GG_{4}$; according to \rref{impifi} we have
$\sup_{k \in \Zt_0} \Gg_{4}(k) \leqs
\sup_{k \in \Zt_0} \GG_{4}(k) \leqs \Big(\sup_{k \in \Zt_0} \Gg_{4}(k) \Big) +
\dG_{4}$, and the numerical results \rref{dg4} \rref{summag4} give
\beq 56.628 < \sup_{k \in \Zt_0} \GG_{4}(k) < 57.892~. \label{sugit4} \feq
\textbf{The upper bound $\boma{G^{+}_{4}}$.}
According to the definition \rref{desup}, we have
\beq G^{+}_{4} = {1 \over (2 \pi)^{3/2}}
\sqrt{\sup_{k \in \Zt_0} \GG_{4}(k)} \quad \mbox{(or any upper approximant for this)}~. \feq
Due to \rref{sugit4}, we can take $G^{+}_4 = (2 \pi)^{-3/2} \sqrt{57.892}\,$; rounding up to  three digits we can write
\beq G^{+}_{4} =  0.484  \, , \feq
as reported in \rref{boug}. \salto
\textbf{Some details on the
evaluation of $\boma{\GG_{5}}$ and of its sup.}
Eq. \rref{dedeg} (with the value of $C_5$
in \rref{valc2}) gives
\beq \dG_{5} = 0.067895...~, \label{dg5} \feq
and it remains to evaluate the function $\Gg_{5}$.
\parn
To compute  $\Gg_{5}(k)$, we start from the $k's$ in $\Zt_0$ with $|k| < 2 \ro = 20$.
Using directly the definition \rref{gmnd}
for all such $k$'s, we obtain
\beq \max_{k \in \Zt_0, |k| < 20} \Gg_{5}(k) = \Gg_{5}(2,1,0) = 138.96...~.  \label{c5} \feq
Let us pass to the case $|k| \geqs 20$. Here we use the upper bound in \rref{tkm} with $t=6$,
giving
\parn
\vbox{
\beq \Gg_5(k) \leqs
40.612 + {271.13 \over |k|^2}
+ {1970.7 \over |k|^4} - {43.608 \over |k|^5}+
{1.4210 \times 10^6 \over |k|^6} - {8949.1 \over |k|^7} \label{e5} \feq
$$  - {2.4425 \times 10^{6} \over |k|^9}
+ {1.6428 \times 10^{5} \over |k|^{10}} - {2.9673 \times 10^{10} \over |k|^{11}}
+ {1.2866 \times 10^{5} \over |k|^{12}}
$$
$$ + {5.3524 \times 10^{5} \over |k|^{14}}
+ {7.9455 \times 10^{14} \over |k|^{16}} \leqs 41.325
\qquad \mbox{for $k \in \Zt_0$, $|k| \geqs 20$}~. $$
}
\noindent
(For completeness we mention that the $t=6$ lower bound in
\rref{tkm} and Eq.\rref{limsup} imply
$\inf_{k \in \Zt_0, |k| \geqs 20} \Gg_{5}(k) = \liminf_{k \in \Zt_0,
k \vain \infty} \Gg_{5}(k) = 8.5405...$ and $\limsup_{k \in \Zt_0,
k \vain \infty} \Gg_{5}(k) = 40.611...\,.$)
\parn
The results \rref{c5} \rref{e5} yield
\beq \sup_{k \in \Zt_0} \Gg_{5}(k) = \Gg_{5}(2,1,0) = 138.96...~. \label{summag5} \feq
We now pass to the function $\GG_{5}$; according to \rref{impifi} we have
$\sup_{k \in \Zt_0} \Gg_{5}(k) \leqs
\sup_{k \in \Zt_0} \GG_{5}(k) \leqs \Big(\sup_{k \in \Zt_0} \Gg_{5}(k) \Big) +
\dG_{5}$, and the numerical results \rref{dg5} \rref{summag5} give
\beq 138.96 < \sup_{k \in \Zt_0} \GG_{5}(k) < 139.04~. \label{sugit5} \feq
\textbf{The upper bound $\boma{G^{+}_{5}}$.}
According to the definition \rref{desup}, we have
\beq G^{+}_{5} = {1 \over (2 \pi)^{3/2}}
\sqrt{\sup_{k \in \Zt_0} \GG_{5}(k)} \quad \mbox{(or any upper approximant for this)}~. \feq
Due to \rref{sugit4}, we can take $G^{+}_5 = (2 \pi)^{-3/2} \sqrt{139.04}\,$; rounding up to  three digits we can write
\beq G^{+}_{5} = 0.749   \, , \feq
as reported in \rref{boug}. \salto
\textbf{Some details on the
evaluation of $\boma{\GG_{10}}$ and of its sup.}
Eq. \rref{dedeg} (with the value of $C_{10}$
in \rref{valc2}) gives
\beq \dG_{10} = 1.0366...\times 10^{-7}~, \label{dg10} \feq
and it remains to evaluate the function $\Gg_{10}$.
\parn
To compute  $\Gg_{10}(k)$, we start from the $k's$ in $\Zt_0$ with $|k| < 2 \ro = 20$.
Using directly the definition \rref{gmnd}
for all such $k$'s, we obtain
\beq \max_{k \in \Zt_0, |k| < 20} \Gg_{10}(k) = \Gg_{10}(2,1,0) = 1.4143...\times 10^4~.  \label{c10} \feq
Let us pass to the case $|k| \geqs 20$. Here we use the upper bound in \rref{tkm} with $t=6$,
giving
\parn
\vbox{
\beq \Gg_{10}(k) \leqs
137.62 + {3125.7 \over |k|^2}
+ {3.2133 \times 10^4 \over |k|^4} + {5.9819 \times 10^7 \over |k|^6}
- {9.2610 \over |k|^{10}}
\label{e10} \feq
$$ - {78.735 \over |k|^{12}}
- {1.1360 \times 10^{4} \over |k|^{14}}
- {1.0781 \times 10^{9} \over |k|^{16}} + {1.6428 \times 10^{5} \over |k|^{20}}
+ {4.9586 \times 10^{8} \over |k|^{22}}
$$
$$ + {6.8396. \times 10^{11} \over |k|^{24}} + {5.0800 \times 10^{17} \over |k|^{26}} \leqs 146.57 \quad
\mbox{for $k \in \Zt_0$, $|k| \geqs 20$}~. $$
}
\noindent
(For completeness we mention that the $t=6$ lower bound in
\rref{tkm} and Eq.  \rref{limsup} imply
$\inf_{k \in \Zt_0, |k| \geqs 20} \Gg_{10}(k)$
$= \liminf_{k \in \Zt_0, k \vain \infty} \Gg_{10}(k) =
4.4157...$ and $\limsup_{k \in \Zt_0,
k \vain \infty}$ $\Gg_{10}(k) = 137.61...\,.$) \parn
The results \rref{c10} \rref{e10} yield
\beq \sup_{k \in \Zt_0} \Gg_{10}(k) = \Gg_{10}(2,1,0) = 1.4143...\times 10^4~. \label{summag10} \feq
We now pass to the function $\GG_{3}$; according to \rref{impifi} we have
$\sup_{k \in \Zt_0} \Gg_{10}(k) \leqs
\sup_{k \in \Zt_0} \GG_{10}(k) \leqs \Big(\sup_{k \in \Zt_0} \Gg_{10}(k) \Big) +
\dG_{10}$, and the numerical results \rref{dg10} \rref{summag10} give
({\footnote{In the MATHEMATICA output for $\Gg_{10}(2,1,0)$,
$1.4143$ is followed by a digit different from $9$; so, the digits $1.4143$ do
not change when $\dG_{10}$ is added to this output.}})
\beq \sup_{k \in \Zt_0} \GG_{10}(k) = 1.4143 ...
\times 10^4~. \label{sugit10} \feq
\textbf{The upper bound $\boma{G^{+}_{10}}$.}
According to the definition \rref{desup}, we have
\beq G^{+}_{10} = {1 \over (2 \pi)^{3/2}}
\sqrt{\sup_{k \in \Zt_0} \GG_{10}(k)} \quad \mbox{(or any upper approximant for this)}~. \feq
Using \rref{sugit10}, and rounding up to  three digits the final result, we can write
\beq G^{+}_{10} =  7.56  \, , \feq
as reported in \rref{boug}.
\vfill \eject \noindent
\section{Appendix. The lower bounds $\boma{G^{-}_{n}}$, for
$\boma{d=3}$ and $\boma{n=3,4,5,10}$}
\label{appe345lowg}
\vskip 1cm
\noindent
Let $n \in (5/2, + \infty)$; according to Proposition \ref{prolowg}, for all nonzero families
$(v_k)_{k \in V}$, $(w_k)_{k \in W}$ in the space $\S$ of \rref{esseu}, we have
the lower bound \rref{gnlow}
\parn
\vbox{
$$ G^{-}_n := {1 \over (2 \pi)^{3/2}} \, {|P_n((v_k), (w_k))| \over N_n((v_k)) N^2_n((w_k))}~
\mbox{(or any lower approximant for this)}, $$
$$ N_n((v_k)) := \Big(\sum_{k \in V} |k|^{2 n} |v_k|^2 \Big)^{1/2}\!\!\!\!\!\!, \qquad
N_n((w_k)) := \Big(\sum_{k \in V} |k|^{2 n} |w_k|^2 \Big)^{1/2}\!\!\!\!\!\!, $$
$$ P_n((v_k), (w_k)) := - i \sum_{h \in V, \ell \in W, h + \ell \in W}
|h + \ell|^{2 n} (\overline{v_h} \sc \ell) (\overline{w_{\ell}} \sc w_{h + \ell})~. $$
}
\noindent
Let us consider the choices
\beq V := \{ \pm (1,0,0) \}~, \qquad v_{\pm (1,0,0)} := (0,P \pm i Q,0) \quad (P, Q \in \reali)~; \feq
\parn
\vbox{
\beq W := \{ \pm (0, 1, 0), \pm (1, 1, 0), \pm (1, -1, 0), \pm (2, 1, 0), \pm (2, -1, 0) \}~; \feq
$$ w_{\pm \k} := (0,0, X_\k \pm i Y_\k) \quad (X_\k, Y_\k \in \reali) $$
$$ ~\mbox{for $\k=(0, 1, 0), (1, 1, 0), (1, -1, 0), (2, 1, 0), (2, -1, 0)$} $$
}
\noindent
(with $(P,Q) \neq 0$ and $(X_\k, Y_\k)_{\k=(0,1,0),...,(2,-1,0)} \neq 0$).
For any $n$, the expressions of
$N_n((v_k))$, $N_n((w_k))$ and $P_n((v_k), (w_k))$
can be computed from the above definitions.
One gets
\beq  N^2_n((v_k)) = 2(P^2 + Q^2)~, \label{esp1} \feq
$$ N^2_n((w_k)) = 2 (X^2_{(0,1,0)} + Y^2_{(0,1,0)}) + 2^{n+1} \hspace{-0.3cm} \sum_{\k=(1,\pm 1,0)} (X^2_\k + Y^2_\k)
+ 2 \times 5^{n} \hspace{-0.3cm} \sum_{\k=(2,\pm 1,0)} (X^2_\k + Y^2_\k)~; $$
\parn
\vbox{
\beq P_n((v_k), (w_k)) = \label{esp2} \feq
$$
2 \Big( - Q X_{(0, 1, 0)} X_{(1, -1, 0)} + Q X_{(0, 1, 0)}
  X_{(1, 1, 0)} + P X_{(1, -1, 0)} Y_{(0, 1, 0)} +
 P X_{(1, 1, 0)} Y_{(0, 1, 0)} $$
$$  + P X_{(0, 1, 0)} Y_{(1, -1, 0)} + Q Y_{(0, 1, 0)} Y_{(1, -1, 0)} -
 P X_{(0, 1, 0)} Y_{(1, 1, 0)} + Q Y_{(0, 1, 0)}
  Y_{(1, 1, 0)} \Big)  $$
$$ + 2^{n+1} \Big(
Q X_{(0, 1, 0)} X_{(1, -1, 0)} -
 Q X_{(0, 1, 0)} X_{(1, 1, 0)} - Q X_{(1, -1, 0)}
  X_{(2, -1, 0)} + Q X_{(1, 1, 0)} X_{(2, 1, 0)} $$
$$ -  P X_{(1, -1, 0)} Y_{(0, 1, 0)} - P X_{(1, 1, 0)} Y_{(0, 1, 0)}
- P X_{(0, 1, 0)} Y_{(1, -1, 0)} -
 P X_{(2, -1, 0)} Y_{(1, -1, 0)} $$
$$ -  Q Y_{(0, 1, 0)} Y_{(1, -1, 0)} + P X_{(0, 1, 0)} Y_{(1, 1, 0)}
+ P X_{(2, 1, 0)} Y_{(1, 1, 0)} -
 Q Y_{(0, 1, 0)} Y_{(1, 1, 0)} $$
$$ + P X_{(1, -1, 0)} Y_{(2, -1, 0)} - Q Y_{(1, -1, 0)} Y_{(2, -1, 0)}
- P X_{(1, 1, 0)} Y_{(2, 1, 0)} + Q Y_{(1, 1, 0)} Y_{(2, 1, 0)}
\Big) $$
$$ + 2 \times 5^n \Big(
Q X_{(1, -1, 0)} X_{(2, -1, 0)} - Q X_{(1, 1, 0)}
  X_{(2, 1, 0)} + P X_{(2, -1, 0)} Y_{(1, -1, 0)} -
 P X_{(2, 1, 0)} Y_{(1, 1, 0)} $$
$$ - P X_{(1, -1, 0)} Y_{(2, -1, 0)} + Q Y_{(1, -1, 0)} Y_{(2, -1, 0)} +
 P X_{(1, 1, 0)} Y_{(2, 1, 0)} - Q Y_{(1, 1, 0)}
  Y_{(2, 1, 0)}
\Big)~. $$
}
\noindent
For any $n$, inserting the expressions \rref{esp1} \rref{esp2} into Eq. \rref{gnlow}
we get a lower bound $G^{-}_n$ depending on the real variables $P, Q, X_\k, Y_\k$. Of course,
to get the best lower bound of this type one should choose $P,Q,X_\k, Y_\k$ so as to maximize
the ratio ${|P_n((v_k), (w_k))|/N_n((v_k)) N^2_n((w_k))}$ in the right hand side
of \rref{gnlow}. \parn
A search of the maximum has been done for $n=3,4,5,10$, using the maximization
algorithms of MATHEMATICA. The program suggests that the maxima
should be attained close to the points $(P,Q,X_\k, Y_\k)$
reported below. It is not granted that such values actually produce
the wanted maxima; in any case, the numbers obtained from \rref{gnlow}
with these choices of $P,Q, X_\k, Y_\k$ are lower bounds on $G_n$, and are
the best derivable by the above algorithms. \parn
The values provided by MATHEMATICA are as follows:
\parn
\vbox{
\beq n = 3: \qquad P = 1, Q = -7.0796..., \feq
$$ X_{(0, 1, 0)} = 1, Y_{(0, 1, 0)} = -5.8246...,
X_{(1, -1, 0)} = -0.063853..., Y_{(1, -1, 0)} = -2.1489..., $$
$$ X_{(1, 1, 0)} = 0.65657..., Y_{(1, 1, 0)} = -2.0472...,
X_{(2, -1, 0)} = -0.043617..., Y_{(2, -1, 0)} = 0.39270..., $$
$$ X_{(2, 1, 0)} = 0.17210..., Y_{(2, 1, 0)} = -0.35566...~; $$
}
\noindent
\parn
\vbox{
\beq n = 4: \qquad P = 1, Q = -7.0768..., \feq
$$  X_{(0, 1, 0)} = 1, Y_{(0, 1, 0)} = -2.7437...,
X_{(1, -1, 0)} = -0.16319..., Y_{(1, -1, 0)} = -0.76896..., $$
$$ X_{(1, 1, 0)} = 0.36987..., Y_{(1, 1, 0)} = -0.69363...,
X_{(2, -1, 0)} = 0.0065160..., Y_{(2, -1, 0)} = 0.094627..., $$
$$  X_{(2, 1, 0)} = 0.055900..., Y_{(2, 1, 0)} = -0.076628...~; $$
}
\noindent
\parn
\vbox{
\beq n = 5: \qquad P = 1, Q = -7.0768..., \feq
$$  X_{(0, 1, 0)} = 1, Y_{(0, 1, 0)} = -2.7618...,
X_{(1, -1, 0)} = -0.12151..., Y_{(1, -1, 0)} = -0.57858..., $$
$$ X_{(1, 1, 0)} = 0.27707..., Y_{(1, 1, 0)} = -0.52225...,
X_{(2, -1, 0)} = 0.0031227..., Y_{(2, -1, 0)} = 0.046786..., $$
$$  X_{(2, 1, 0)} = 0.027554..., Y_{(2, 1, 0)} = -0.037939...~; $$
}
\noindent
\parn
\vbox{
\beq n = 10: \qquad P = 1, Q = -7.0769..., \feq
$$  X_{(0, 1, 0)} = 1, Y_{(0, 1, 0)} =  -2.8038...,
X_{(1, -1, 0)} = -0.031443..., Y_{(1, -1, 0)} = -0.15337..., $$
$$ X_{(1, 1, 0)} = 0.072707..., Y_{(1, 1, 0)} = -0.13865...,
X_{(2, -1, 0)} = 8.9903 \times 10^{-5}..., Y_{(2, -1, 0)} = 0.0014520..., $$
$$  X_{(2, 1, 0)} = 8.4924 \times 10^{-4}..., Y_{(2, 1, 0)} = - 0.0011812...~. $$
}
\noindent
(Note that the ratio $|P_n((v_k), (w_k))|/N_n((v_k)) N^2_n((w_k))$
is invariant under any rescaling $(v_k) \mapsto (\lambda v_k)$,
$(w_k) \mapsto (\mu w_k)$, with $\lambda, \mu \in \reali \setminus \{0\}$;
the normalizations for $P$ and $X_{(0,1,0)}$ adopted above arise from
the possibility of such rescalings.) \parn
With the above choices of $P,Q, X_\k, Y_\k$ (i.e., of $(v_k)$ and
$(w_k)$), one has
\beq G^{-}_n ~=
\left\{ \barray{cc} 0.11433... & \mbox{for $n=3$,} \\ 0.18128... & \mbox{for $n=4$,} \\
0.28013... & \mbox{for $n=5$,} \\ 2.4155... & \mbox{for $n=10$.} \farray \right. \feq
Rounding down to three digits the above numbers,
we obtain the results in \rref{boug}.
\vskip 0.7cm \noindent
\textbf{Acknowledgments.} This work was partly supported by INdAM and by MIUR, PRIN 2008
Research Project "Geometrical methods in the theory of nonlinear waves and applications"

\end{document}